\documentclass[twocolumn,twoside]{IEEEtran}
\usepackage{amsmath,amssymb,amsthm,epsfig,color,subfigure,empheq,graphicx,graphics,pgfplots}
\usepackage{enumerate,url,algpseudocode,algorithm,wasysym,epstopdf,balance,enumitem}
\usepackage{multirow} 

\newtheorem{proposition}{Proposition}

\newtheorem{remark}{Remark}

\newcommand \bzero{\mathbf{0}}
\newcommand \bone{\mathbf{1}}
\newcommand \ba{\mathbf{a}}
\newcommand \bb{\mathbf{b}}

\newcommand \bd{\mathbf{d}}
\newcommand \be{\mathbf{e}}
\newcommand \bef{\mathbf{f}} 

\newcommand \bp{\mathbf{p}}
\newcommand \bq{\mathbf{q}}

\newcommand \bu{\mathbf{u}}
\newcommand \bv{\mathbf{v}}
\newcommand \bw{\mathbf{w}}
\newcommand \bx{\mathbf{x}}
\newcommand \by{\mathbf{y}}
\newcommand \bz{\mathbf{z}}

\newcommand \bI{\mathbf{I}}

\newcommand \bK{\mathbf{K}}

\newcommand \bM{\mathbf{M}}

\newcommand \bR{\mathbf{R}}

\newcommand \bX{\mathbf{X}}

\newcommand \bZ{\mathbf{Z}}

\newcommand \bgamma{\boldsymbol{\gamma}}

\newcommand \bmu{\boldsymbol{\mu}}

\newcommand \brho{\boldsymbol{\rho}}
\newcommand \bsigma{\boldsymbol{\sigma}}

\newcommand \bphi{\boldsymbol{\phi}}
\newcommand \bchi{\boldsymbol{\chi}}

\newcommand \bPhi{\mathbf{\Phi}}

\newcommand \mcH{\mathcal{H}}

\newcommand \mcK{\mathcal{K}}

\newcommand \mcN{\mathcal{N}}
\newcommand \mcO{\mathcal{O}}
\newcommand \mcP{\mathcal{P}}
\newcommand \mcQ{\mathcal{Q}}

\newcommand \mcS{\mathcal{S}}

\newcommand \mcZ{\mathcal{Z}}

\newcommand \tbq{\tilde{\mathbf{q}}}

\newcommand \ublambda{\underline{\boldsymbol{\lambda}}}
\newcommand \oblambda{\overline{\boldsymbol{\lambda}}}
\newcommand \obmu{\overline{\boldsymbol{\mu}}}
\newcommand \ubmu{\underline{\boldsymbol{\mu}}}

\begin{document}
	\title{Designing Reactive Power Control Rules\\
	for Smart Inverters using Support Vector Machines}
	
\author{
		Mana Jalali,~\IEEEmembership{Student Member,~IEEE},
		Vassilis Kekatos,~\IEEEmembership{Senior Member,~IEEE}, 		Nikolaos Gatsis,~\IEEEmembership{Member,~IEEE}, and 
		Deepjyoti Deka,~\IEEEmembership{Member,~IEEE},
		
\thanks{Manuscript received March 1, 2019; revised June 4, 2019, and August 15, 2019; accepted September 16, 2019. Date of publication DATE; date of current version DATE. Paper no. TSG.00318.2019. This work was supported in part by the NSF-CAREER grant 1751085.}
		
\thanks{M. Jalali and V. Kekatos are with the Bradley Dept. of ECE, Virginia Tech, Blacksburg, VA 24061, USA. N. Gatsis is with the ECE Dept., University of Texas at San Antonio, San Antonio, TX 78249, USA. D. Deka is with the Theoretical Division at Los Alamos National Laboratory, NM 87545, USA. Emails: \{manaj2,kekatos\}@vt.edu; nikolaos.gatsis@utsa.edu; deepjyoti@lanl.gov.}	
\thanks{Color versions of one or more of the figures is this paper are available online at {http://ieeexplore.ieee.org}.}
\thanks{Digital Object Identifier XXXXXX}
}	
	
	\markboth{IEEE TRANSACTIONS ON SMART GRID (to appear)}{Jalali, Kekatos, Gatsis, and Deka: Designing Reactive Power Control Rules for Smart Inverters using Support Vector Machines}
	
	\maketitle
	
	\begin{abstract}
Smart inverters have been advocated as a fast-responding mechanism for voltage regulation in distribution grids. Nevertheless, optimal inverter coordination can be computationally demanding, and preset local control rules are known to be subpar. Leveraging tools from machine learning, the design of customized inverter control rules is posed here as a multi-task learning problem. Each inverter control rule is modeled as a possibly nonlinear function of local and/or remote control inputs. Given the electric coupling, the function outputs interact to yield the feeder voltage profile. Using an approximate grid model, inverter rules are designed jointly to minimize a voltage deviation objective based on anticipated load and solar generation scenarios. Each control rule is described by a set of coefficients, one for each training scenario. To reduce the communication overhead between the grid operator and the inverters, we devise a voltage regulation objective that is shown to promote parsimonious descriptions for inverter control rules. Numerical tests using real-world data on a benchmark feeder demonstrate the advantages of the novel nonlinear rules and explore the trade-off between voltage regulation and sparsity in rule descriptions. 
	\end{abstract}
	
	
	\allowdisplaybreaks
	
	\section{Introduction}\label{sec:intro}
	Several electric utilities in the US currently experience issues while integrating residential- and commercial-scale solar generation. A solar farm connected at the end of a long rural feeder can incur voltage excursions along the feeder, while frequent power flow reversals strain the apparent power capabilities of substation transformers~\cite{Turitsyn11}. Solar generation from residential photovoltaics (PVs) can fluctuate by up to 15\% of their rating within one-minute intervals~\cite{pecandata}. Utility-owned voltage control equipment, such as load-tap-changing transformers, capacitor banks, and step-voltage regulators, involves discrete control actions, and its lifespan is related to the number of switching operations~\cite{tse2015jabr}. Regulating voltage under increasing renewable generation may require more frequent switching and further installations, thus critically challenging reactive power control in distribution grids. 
	
	On the other hand, PVs are interfaced by inverters featuring advanced communication, metering, and control functionalities. Using inverters for reactive power control has been advocated as a fast-responding solution~\cite{Turitsyn11}. The amended IEEE 1547 standard allows inverters to be operating at non-unit power factors~\cite{IEEE1547}. Nonetheless, coordinating in real-time hundreds of inverters distributed over a feeder is a formidable task. In a typical setup, the values of instantaneous loads and solar generation are communicated to a utility controller; the controller minimizes ohmic losses subject to voltage regulation constraints; and the computed setpoints are sent back to inverters. The problem of finding the optimal reactive injection setpoints for inverters is an instance of the optimal power flow (OPF) task, which is non-convex in general. Different convex relaxations have been proposed; see \cite{Low14} for a survey. The uncertainty in loads and solar generation over the next control period is usually accounted for through stochastic and robust formulations~\cite{ergodic}, \cite{YGG13}. To reduce complexity, approximate grid models have also been employed \cite{FCL}, \cite{Robbins}; though heavy two-way utility-inverter communication is still needed. 
	
	Alternatively, decentralized solvers where inverters decide their setpoints upon communicating with neighboring inverters have been devised~\cite{tse2014emiliano}, \cite{QiuyuPeng}, \cite{BaGa-TSG2017}. On the other extreme, localized schemes suggest having inverters implementing Volt-VAR and/or Watt-VAR curves given only local measurements~\cite{Turitsyn11}. Although such rules have been analytically shown to be stable and fast-converging, their equilibria unfortunately do not coincide with the sought OPF minimizers~\cite{ZDGT13}, \cite{FCL}, \cite{LQD14}, \cite{VKZG16}. In fact, there exist cases where local rules perform worse than the no-reactive support option~\cite{Jabr18}. 
	
The previous literature review indicates that centralized schemes incur high computational complexity; decentralized solvers require multiple communication exchanges among inverters; and local schemes have no performance guarantees. As a middle-ground solution, inverter setpoints can be designed in a quasi-static fashion via control rules. A rule expresses each setpoint as an affine function of given inputs, such as generation, load, or voltage. Albeit the related weights are optimized periodically in a centralized fashion, control rules are applied in real time. Controlling inverters via affine rules has been accomplished using chance-constrained~\cite{Ayyagari17}; robust \cite{Jabr18}, \cite{LinBitar18}; and closed-loop formulations~\cite{Baker18}. Optimal rules however are not necessarily linear: If an apparent power constraint becomes active, reactive injections can become nonlinear functions of solar generation. To capture this nonlinearity, recent approaches engage learning models which are \emph{trained to optimize}: Given pairs of grid conditions (load and solar generation) and their optimal inverter dispatches computed, the aforesaid approaches learn dispatch rules using linear or kernel-based regression~\cite{Dobbe19}, \cite{Kara18}.

This work combines machine learning tools with physical grid models, and advocate a kernel-based approach for designing inverter control rules. The contribution is on two fronts: First, the design of inverter control rules is posed as a multi-task learning problem. Each inverter rule is modeled as a nonlinear function of control inputs. Rules are coupled through the electric grid to yield a system voltage profile. Using an approximate grid model, inverter rules are learned jointly so that they minimize a voltage regulation cost using anticipated load and solar generation scenarios. Each rule is described by a set of coefficients, one for each scenario. As a second contribution, we engineer the voltage regulation objective, so that the optimal rules are described by a few scenario coefficients. Such parsimonious representation of inverter rules saves communications. Numerical tests on a benchmark feeder showcase the advantages of nonlinear rules and explore the trade-off between voltage regulation and sparse rules. 
	
Regarding notation, lower- (upper-) case boldface letters denote column vectors (matrices), while calligraphic symbols are reserved for sets. Symbol $^{\top}$ stands for transposition and $\|\bx\|_2$ denotes the $\ell_2$-norm of $\bx$.
	
\section{Reactive Power Control}\label{sec:model}
This section formulates the task of voltage regulation using inverters. Consider a distribution grid having $N+1$ buses served by the substation indexed by $n=0$. Let $v_n$ denote the voltage magnitude, and $p_n+jq_n$ the complex power injection at bus $n$. The active injection $p_n$ is decomposed into $p_n=p_n^g-p_n^c$, where $p_n^g$ is the solar generation and $p_n^c$ the inelastic load at bus $n$. Reactive injections can be similarly expressed as $q_n=q_n^g-q_n^c$. Collect injections in $N$-length vectors:
		\begin{equation}\label{eq:pq}
		\bp=\bp^g-\bp^c~~\textrm{and}~~
		\bq=\bq^g-\bq^c.
		\end{equation}
	
	The reactive power injected by inverter $n$ is constrained as
	\begin{equation}\label{eq:pv}
	|q_n^g|\leq \bar{q}_n^g:=\sqrt{(\bar{s}_n^g)^2 - (p_n^g)^2}
	\end{equation}
	where $\bar{s}_n^g$ is the apparent power limit for inverter $n$; see~\cite{Turitsyn11}.
	
	Given loads $(\bp^c,\bq^c)$ and solar generation $\bp^g$, voltage regulation aims at optimally setting $\bq^g$ such that voltage deviations are kept minimal. To formally describe this task, one has to deal with the nonlinear power flow equations relating voltages to power injections. Trading modeling accuracy for computational tractability, we resort to the linearized model~\cite{BoDo15}
	\begin{equation}\label{eq:linear}
	\bv \simeq \bR\bp + \bX\bq + v_0\bone
	\end{equation}
	where $\bv:=[v_1~\ldots~v_N]^\top$ and matrices $(\bR,\bX)$ depend on the feeder. Model \eqref{eq:linear} can be derived by linearizing the power flow equations around the flat voltage profile. In fact, the linearization can be performed at any system state $\bv_0$, yet matrices $\bR$ and $\bX$ would then depend on the state~$\bv_0$; see~\cite{Baker18}. From \eqref{eq:pq} and \eqref{eq:linear}, the vector of voltage deviations from its nominal value can be approximated as
	\begin{equation}\label{eq:dv}
	\bv-v_0\bone = \bX\bq^g +\by
	\end{equation}
	where $\by:=\bR(\bp^g-\bp^c) - \bX\bq^c$ and $\bone$ is a vector of all ones. 

The goal here is to design the inverter injections $\bq^g$ so that bus voltage magnitudes remain within regulation limits. The ANSI-C.84.1 standard dictates that service (load) voltages should remain within $\pm 5\%$~per unit (pu). However, our grid model of \eqref{eq:dv} stops at the level of distribution transformers. A distribution (pole or pad-mounted) transformer may be serving several residential customers. Each customer is typically connected to the distribution transformer through a triplex cable, which incurs a voltage drop between the transformer and the service voltage: Suppose a customer is connected to a 50~kVA, 7200-240/120~V center-tapped transformer via a 1/0~AA 100-ft triplex cable. The customer runs a constant-current load of 10~kVA at the nominal voltage of 120~V with 0.9 lagging power factor. If load currents are equally distributed among the three supplies (two 120~V and one 240~V), the service voltage drops by $1.5\%$~pu. If the load is distributed among supplies non-uniformly, the service voltage can drop by even $3.5\%$~pu. Due to this, the current practice is to maintain voltages at distribution transformers within $\pm 3\%$~pu, to ensure that service voltages remain within $\pm 5\%$~pu; see exercises of~\cite{Kersting}.
		
Given loads, solar generation, and grid parameters, the goal is to decide $\bq^g$ to regulate voltage while satisfying the apparent power constraints of \eqref{eq:pv}. The setpoints for reactive power injections from inverters can be found as the minimizer
	\begin{align}\label{eq:vr}
	\tbq^g:=\arg\min_{\bq^g\in \mcQ}&~\Delta(\bq^g;\by).
	\end{align}
	The set $\mcQ\subseteq \mathbb{R}^N$ captures the constraints in \eqref{eq:pv} for all $n$; and $\Delta(\bq^g;\by)$ is a voltage regulation objective. A typical choice for $\Delta$ is the sum of squared voltage deviations~\cite{ZDGT13}, \cite{Jabr18}, \cite{Baker18}
	\begin{equation}\label{eq:Ds}
	\Delta_s(\bq^g;\by):=\sum_{n=1}^N\left(v_n-v_0\right)^2=\|\bX\bq^g+\by\|_2^2.
	\end{equation}
Alternatively, the utility may want to maintain voltages within the range of $(1\pm\epsilon)v_0$ for say $\epsilon=0.03$. Then, a pertinent objective is~\cite{LQD14}
	\begin{equation}\label{eq:De}
	\Delta_\epsilon(\bq^g;\by):=\sum_{n=1}^N \left[v_n-v_0\right]_\epsilon=\sum_{n=1}^N \left[\be_n^\top\left(\bX \bq^g +\by\right)\right]_\epsilon
	\end{equation}
	where $\be_n$ is the $n$-th canonical vector of length $N$, and the operator $[\cdot]_\epsilon$ is defined as
	\begin{equation}\label{eq:epsilon}
	[x]_\epsilon:=\left\{\begin{array}{ll}
	0 &,~|x|\leq \epsilon\\
	|x|-\epsilon&,~\text{otherwise}
	\end{array}\right..
	\end{equation}
Function $\Delta_\epsilon$ returns zero when all voltages are within limits. Otherwise, it increases linearly with voltage excursions; see~\cite{LQD14} for distributed solvers of \eqref{eq:vr} with $\Delta=\Delta_\epsilon$.
	
	It is worth noticing that $\bX$ depends only on the network and the linearization point, whereas the set $\mcQ$ and vector $\by$ depend on the variable loads and solar generation, collectively denoted as $\bchi:=[(\bp^c)^\top~(\bq^c)^\top~(\bp^g)^\top]^\top$. 
	
	Ideally, the reactive control process entails three steps:
	\renewcommand{\labelenumi}{\emph{S\arabic{enumi})}}
	\begin{enumerate}
		\item Each bus communicates its $(p_n^g,p_n^c,q_n^c)$ to the operator.
		\item The operator solves \eqref{eq:vr} knowing the current $\bchi$.
		\item The operator sends the optimal setpoints $\tbq^g$ to inverters.
	\end{enumerate}
	
	Under variable solar generation, the process has to be repeated on a per-minute basis. Observe that \emph{S1)} establishes $N$ inverter-to-utility communication links, and \emph{S3)} requires another $N$ utility-to-inverter links. Running this process for multiple feeders can become a computationally and communication-wise challenging task.
	
	To adaptively adjust inverter setpoints based on $\bchi_t$, affine control rules in the form of $\bq^g(\bchi_t)$ have been suggested in \cite{Jabr18}, \cite{Ayyagari17}, \cite{LinBitar18}. Based on these rules, the reactive injection of inverter $n$ is expressed as an affine function over a subvector of $\bchi_t$. The premise is to design the rule in a quasi-stationary fashion, but apply it in real-time. We extend linear to \emph{nonlinear} control rules enjoying varying cyber requirements after briefly reviewing the toolbox of kernel-based learning.
	
\section{Preliminaries on Kernel-based Learning}\label{sec:background}
Given pairs $\{(z_s,y_s)\}_{s=1}^S$ of features $z_s$ belonging to a measurable space $\mcZ$ and target values $y_s\in \mathbb{R}$, kernel-based learning aims at finding a function or mapping $f:\mcZ\rightarrow\mathbb{R}$. From all possible options of arbitrarily complex functions, one needs to select a specific family where $f$ belongs. Kernel-based learning postulates that $f$ lies in the function space~\cite{Hastie}
	\begin{equation}\label{eq:family}
	\mcH_\mcK:=\left\{f(z)=\sum_{s=1}^{\infty} K(z,z_s) a_s,~a_s\in\mathbb{R}\right\}.
	\end{equation}
This is the space of functions that can be expressed as linear combinations of a given kernel (basis) function $K:\mathcal{Z}\times \mathcal{Z}\rightarrow \mathbb{R}$ evaluated at arbitrary points $z_s$. When $K(\cdot,\cdot)$ is a symmetric positive definite function, then $\mathcal{H}_{\mathcal{K}}$ becomes a reproducing kernel Hilbert space (RKHS) whose members have finite norm $\|f\|_{\mathcal{K}}^2:= \sum_{s=1}^{\infty} \sum_{s'=1}^{\infty} K(z_s,z_{s'}) a_s a_{s'}$; see~\cite{Aber09}. Some options for the kernel function $K$ are provided under \emph{Examples 1--2} in Section~\ref{subsec:rules2}. 
	
Learning $f$ from data $\{(z_s,y_s)\}_{s=1}^S$ can be formulated as the regularization task~\cite{Hastie}, \cite{BaGia13}
	\begin{equation}\label{eq:fnreg}
	\min_{f\in\mathcal{H}_{\mathcal{K}},b}~ \frac{1}{S}\sum_{s=1}^S L\left(f(z_s),b;y_s\right) + \mu \|f\|_{\mathcal{K}}
	\end{equation}
where $b$ is an intercept term. When it comes to regression, typical choices for the data-fitting loss $L$ include the least-squares (LS) fit $\left(y_s-f(z_s)-b\right)^2$, or the $\epsilon$-insensitive loss $\left[y_s-f(z_s)-b\right]_\epsilon$. The second term in \eqref{eq:fnreg} ensures $f\in\mathcal{H}_{\mathcal{K}}$ and facilitates generalization over unseen data~\cite{Aber09}. Parameter $\mu>0$ balances fitting versus generalization, and is tuned via cross-validation: \emph{i)} problem \eqref{eq:fnreg} is solved for a specific $\mu$ using $4/5$ of the data; \emph{ii)} the learned function is validated on the unused $1/5$ of the data; \emph{iii)} the process is repeated $5$ times to calculate the average fitting error for this $\mu$; and \emph{iv)} the $\mu$ attaining the best fit is selected; see~\cite{Hastie} for details.
	
The advantage of confining $f$ to lie in the RKHS $\mcH_\mcK$ is that the functional optimization of \eqref{eq:fnreg} can be equivalently posed as an minimization problem over a finite-dimensional vector: The celebrated Representer's Theorem asserts that the solution to \eqref{eq:fnreg} admits the form~\cite{Hastie}
	\begin{equation}\label{eq:rt}
	f(z)=\sum_{s=1}^S K(z,z_s)a_s.
	\end{equation}
	In other words, the minimizer of \eqref{eq:fnreg} is described only by $S$ rather than infinitely many $a_s$'s. 
	Based on \eqref{eq:rt}, evaluating $f(z)$ at the given data provides $\bef = \bK \ba$, where $\bef:=[f(z_1)~\ldots~f(z_S)]^\top$; matrix $\bK\in \mathbb{S}_{++}^{S}$ is the \emph{kernel matrix} with entries $[\bK]_{s,s'}:=K(z_s,z_{s'})$; and $\ba:=[a_1~\ldots~a_S]^\top$.
	
	From properties of the RKHS's, it holds that $\|f\|_\mcK^2=\ba^\top\bK\ba$; see \cite{Aber09}. For regression under an LS loss, the functional minimization in \eqref{eq:fnreg} becomes the vector optimization
	\begin{equation}\label{eq:vecreg}
	\min_{\mathbf{a},b}~\frac{1}{S}\|\mathbf{y}-\mathbf{K}\mathbf{a}-b\bone\|_2^2 + \mu\|\bK^{1/2}\ba\|_2
	\end{equation}
where $\bK^{1/2}$ is the square root of $\bK$ and $\mathbf{y}:=[y_1 ~\cdots~y_S]^\top$. 
	
It is worth stressing that \eqref{eq:rt} applies not only to the given data $\{z_s\}_{s=1}^S$, but any $z_{s'}\in\mcZ$. Evaluating $f(z)$ requires knowing the $(\ba,b)$ minimizing \eqref{eq:vecreg}, and being able to evaluate the kernel $K(z,z_s)$ for $s=1,\ldots,S$. We next use kernel-based learning to develop nonlinear inverter control rules.
	
	\section{Kernel-based Control Policies}\label{subsec:rules}
	The reactive injection by inverter $n$ is modeled by the rule
	\begin{equation}\label{eq:qfun}
	q_n^g(\bz_n)=f_n(\bz_n)+b_n
	\end{equation}
	whose ingredients $(f_n,\bz_n,b_n)$ are explained next. 
	
\emph{Control inputs:} Vector $\bz_n\in\mcZ_n\subseteq \mathbb{R}^{M_n}$ is the input to control rule for inverter $n$. This vector may include load, solar generation, and/or line flow measurements collected locally or remotely. For a purely local rule, this input can be selected as
	\begin{equation}\label{eq:input}
	\bz_n:=\left[\bar{q}_{n}^g\quad (p_n^c-p_n^g) \quad q_n^c\right]^\top
	\end{equation}
	where the first entry $\bar{q}_{n}^g$ relates to the apparent power constraint and has been defined in \eqref{eq:pv}. The voltage $v_n$ could also be appended in $\bz_n$; however the stability of the resultant control loop is hard to analyze even when $f_n$ is linear; see e.g.,~\cite{FCL13}, \cite{LQD14}, \cite{VKZG16}, \cite{Baker18}, \cite{Dobbe19}.

Selecting the controller structure, i.e., the content for each $\bz_n$, can affect critically the performance of this control scheme. Ideally, each inverter rule can be fed all uncertain quantities, that is the three numbers in the right-hand side of \eqref{eq:input} across all buses. In that case, the input vectors $\bz_n$ become all equal and of size $3N$. However, this incurs the communication burden of broadcasting $3N$ values in real time. Hybrid setups with $\bz_n$'s carrying a combination of local and remote data can be envisioned. To eliminate the effect of this trade-off between communications and performance, this work assumes that the content of $\bz_n$'s is prespecified. The task of input selection could be possibly pursued along the lines of sparse linear or polynomial regression~\cite{Baker18}, \cite{Dobbe19}, \cite{VKGG11b}; and automatic relevance determination~\cite[Sec.~6.4]{Bishop}.

\emph{Control function:} Selecting the form of $f_n$ is the second design task. To leverage kernel-based learning, the inverter rule $f_n$ is postulated to lie in the RKHS
	\begin{equation}\label{eq:rpcfamily}
	\mcH_{\mcK_n}:=\left\{f_n(\bz_n)=\sum_{s=1}^{\infty} K_n(\bz_{n},\bz_{n,s}) a_{n,s},~a_{n,s}\in\mathbb{R}\right\}
	\end{equation}
 determined by the kernel function $K_n:\mcZ_n\times \mcZ_n\rightarrow \mathbb{R}$. 

Linear rules can be designed by selecting the linear kernel $K_n(\bz_{n,s},\bz_{n,s'})=\bz_{n,s}^\top\bz_{n,s'}$. \emph{Nonlinear} rules can be designed by selecting for example a polynomial kernel $K_n(\bz_{n,s},\bz_{n,s'})=\left(\bz_{n,s}^\top\bz_{n,s'}+\gamma\right)^\beta$ or a Gaussian kernel $K_n(\bz_{n,s},\bz_{n,s'})=\exp\left(-\|\bz_{n,s}-\bz_{n,s'}\|_2^2/\gamma\right)$ with design parameters $\beta>0$ and $\gamma>0$; see~\cite{Hastie}.
	
	\emph{Intercept $b_n\in\mathbb{R}$:} Although it could be incorporated into $f_n$ by augmenting $\bz_n$ with a constant entry of $1$, it is kept separate to avoid its penalization through $\|f\|_{\mcK_n}$~\cite{Hastie}.
	
\subsection{Learning rules from scenario data}\label{subsec:rules2}
The rules of \eqref{eq:qfun} can be learned from scenario data indexed by $s\in \mcS$ with $\mcS:=\{1,\ldots,S\}$. Scenario $s$ consists of the control inputs $\bz_{n,s}$ for $n\in\mcN$, and the associated vector $\by_s:=\bR(\bp^g_s-\bp^c_s) - \bX\bq^c_s$ defined in \eqref{eq:dv}. Evaluating rule $n$ of \eqref{eq:qfun} under scenario $s$ yields the inverter response $q_{n,s}^g:=q_n^g(\bz_{n,s})$. Let us collect the outputs $q_{n,s}^g$ from all inverters into vector $\bq^g_s$. Note that the goal is not to fit $\by_s$ by $\bq^g_s$, but to minimize the voltage deviations $\bX\bq^g_s+\by_s$. The control functions $\{f_n\}_{n=1}^N$ and the intercepts $\{b_n\}_{n=1}^N$ accomplishing this goal can be found via the functional minimization 
	\begin{align}\label{eq:rpcfnreg}
	\min~&~\frac{1}{S}\sum_{s=1}^S \Delta\left(\bq^g_s;\by_s\right) + \mu \sum_{n=1}^N\|f_n\|_{\mcK_n}\\
	\mathrm{over}~&~q_{n,s}^g=f_n(\bz_{n,s})+b_n,~\forall n,s\nonumber\\
	~&~\{f_n\in\mcH_{\mcK_n}\}, \bb:=[b_1~\cdots~b_N]^\top\nonumber\\
	\mathrm{s.to}~&~|q_{n,s}^g|\leq \bar{q}_{n,s}^g,~ \forall n,s\nonumber
	\end{align}
where $\Delta$ is a voltage regulation objective [cf.~\eqref{eq:Ds}--\eqref{eq:De}].

\begin{remark}\label{re:others}
The proposed approach is related to~\cite{Dobbe19}--\cite{Kara18}, where inverter rules are also trained using machine learning. However, the aforementioned works proceed in two steps: They first solve a sequence of OPF problems similar to \eqref{eq:vr} to find the optimal inverter setpoints $\tbq^g$ under different scenarios. Secondly, they learn the mapping between controller inputs $\{\bz_{n,s}\}_{s\in\mcS}$ and optimal setpoints $\{\tilde{q}_{n,s}^g\}$ decided by the OPF problems. During this process, they also select which inputs are more effective to be communicated to inverters. The mapping is learned via linear or kernel-based regression. On the other hand, the approach proposed here consolidates the OPF and the learning steps into a single step: The advantage is that the OPF decisions of \eqref{eq:rpcfnreg} are taken under the explicit practical limitation that $q_{n}^g$ can only be a function of $\bz_{n}$, since inverter $n$ will not have access to the complete grid conditions. To get some intuition, suppose ones designs linear control rules of known input structure using the single-step approach of \eqref{eq:rpcfnreg} with $\mu=0$ and the two-step approach of \cite{Dobbe19}--\cite{Kara18}. The single-step approach yields rules $R_1$, and the two-step approach yields rules $R_2$. Let us evaluate $R_1$ and $R_2$ on the training scenarios. Rules $R_2$ are not necessarily feasible per scenario $s\in\mcS$, whereas rules $R_1$ are. Moreover, rules $R_2$ do not necessarily coincide with the minimizers of \eqref{eq:vr}. For the sake of comparison, let us assume that rules $R_2$ turn out to be feasible per scenario, and hence feasible for \eqref{eq:rpcfnreg}. Being the minimizers of \eqref{eq:rpcfnreg}, rules $R_1$ attain equal or smaller voltage deviation cost compared to $R_2$ over the training data. Numerical tests in Section~\ref{sec:tests} corroborate the advantage of $R_1$ over $R_2$ for $\mu>0$ and during the operational phase as well.
\end{remark}
	
Different from \eqref{eq:fnreg}, the optimization in \eqref{eq:rpcfnreg} entails learning multiple functions (one per inverter). Since inverter injections affect voltages feeder-wise, inverter rules are naturally coupled through $\Delta$ in \eqref{eq:rpcfnreg}. Similar multi-function setups can be found in collaborative filtering or multi-task learning~\cite{Aber09}, \cite{KZG14}. 
	
	Fortunately, Representer's Theorem can be applied successively over $n$ in \eqref{eq:rpcfnreg}. Therefore, each rule $n$ is written as 
\begin{equation}\label{eq:fn}
f_n(\bz_n)=\sum_{s=1}^S K_n(\bz_n,\bz_{n,s})a_{n,s}.
\end{equation}
	Once the coefficients $\{a_{n,s}\}$ have been found, rule $\{f_n\}$ can be evaluated for any $\bz_n$. Similar to \eqref{eq:rt}, evaluating rule $f_n$ over the scenario data  $\{\bz_{n,s}\}_{s=1}^S$ gives 
	\begin{equation}\label{eq:rtn}
	\bef_n = \bK_n \ba_n,\quad \forall n
	\end{equation}
	where $[\bK_n]_{s,s'}=K_n(\bz_{n,s},\bz_{n,s'})$ for $s,s'=1,\ldots,S$, and $\ba_n:=[a_{n,1}~\cdots~a_{n,S}]^\top$. The RKHS norms can be written as
	\begin{equation}\label{eq:norms}
	\|f_n\|_{\mcK_n}=\sqrt{\ba_n^\top\bK_n\ba_n},\quad \forall n.
	\end{equation}
	In this way, the functional minimization in \eqref{eq:rpcfnreg} is cast as a vector minimization over $\{\ba_n\}_{n=1}^N$ and $\bb$. The exact form of this minimization and its properties for different $\Delta$ are discussed later in Section~\ref{sec:svm}. For now, let us clarify how the kernel functions $K_n(\cdot,\cdot)$ effect different rule forms. 
	
	\emph{Example 1: Affine rules.} The linear kernel $K_n(\bz_{n,s},\bz_{n,s'})=\bz_{n,s}^\top\bz_{n,s'}$ yields affine rules. The sought functions can be written as
	\begin{equation}\label{eq:qfunlin}
	f_n(\bz_n)=\bz_n^\top \bw_n,~\quad \forall n.
	\end{equation}
Given scenario data $\bz_{n,s}$ and $\by_s$ for $n\in\mcN$ and $s\in\mcS$, we would like to find $\{\bw_n,b_n\}_n$ through \eqref{eq:rpcfnreg}. Collect the input data for inverter $n$ in the $M_n\times S$ matrix $\bZ_n:=\left[\bz_{n,1}~\cdots~\bz_{n,S}\right]$. According to Representer's Theorem, the optimal $\bw_n$ can be expressed as $\bw_n=\bZ_n\ba_n$ for some $\ba_n$. Evaluating the control rule for any input $\bz_{n,s}$ yields \[q_n(\bz_{n,s})=f_n(\bz_{n,s})+b_n=\bz_{n,s}^\top \bZ_n\ba_n + b_n.\]
	Evaluating the rule at the input data yields \eqref{eq:rtn} with $\bK_n=\bZ_n^\top \bZ_n$. The squared function norm is $\|f_n\|_{\mcK_n}^2=\|\bw_n\|_2^2=\ba_n^\top\bZ_n^\top\bZ_n\ba_n=\ba_n^\top\bK_n\ba_n$.
	
\emph{Example 2: Non-linear rules.} For non-linear rules, transform the input $\bz_{n,s}$ to vector $\bphi_{n,s}:=\phi_n(\bz_{n,s})$ via a non-linear mapping $\phi_n:\mathbb{R}^{M_n}\rightarrow\mathbb{R}^{\Phi_n}$. The entries of $\bphi_{n,s}$ could be for example all the first- and second-order monomials formed by the entries of $\bz_{n,s}$. The dimension $\Phi_n$ of $\bphi_{n,s}$ can be finite (e.g., polynomial kernels) or infinite (Gaussian kernels)~\cite{Bishop}. Then, the control function
	\begin{equation}\label{eq:qfunnlin}
	f_n(\bz_n)=\bphi_n^\top \bw_n
	\end{equation}
	with $\bw_n\in\mathbb{R}^{\Phi_n}$ is \emph{non-linear} in $\bz_n$. The developments of Example~1 carry over to Example 2 by using $\bK_n=\bPhi_n^\top\bPhi_n$ and replacing $\bZ_n$ by $\bPhi_n:=[\bphi_{n,1}~\cdots~\bphi_{n,S}]$. Depending on the mapping $\phi_n$, the vectors $\bphi_{n,s}$ may be of finite or infinite length~\cite{Hastie}. The critical point is that $f_n$ does not depend on $\bphi_{n,s}$'s directly, but only on their inner products $\bphi_{n,s}^\top\bphi_{n,s'}$ for any $s$ and $s'$. These products can be easily calculated through the kernel function as $\bphi_{n,s}^\top\bphi_{n,s'}=K_n(\bz_{n,s},\bz_{n,s'})$; see \cite{Hastie}.
	
	\begin{figure*}[t]
	\centering
	\includegraphics[scale=0.52]{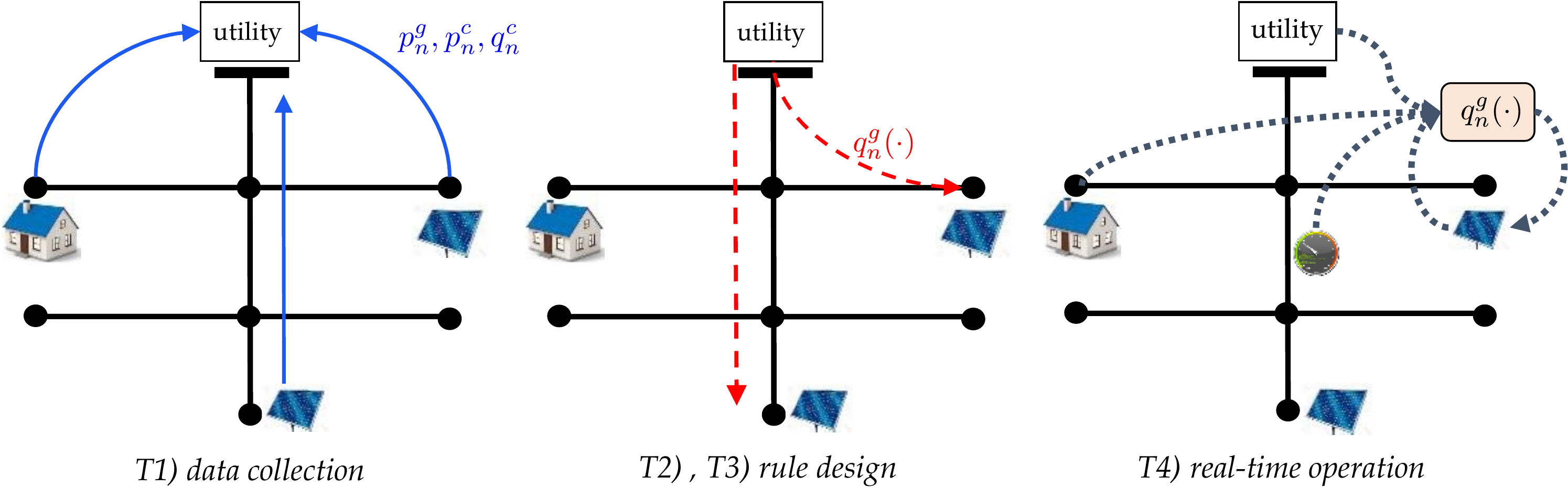}
	\caption{Implementing reactive power control rules. \emph{Left:} Data are collected from buses. \emph{Center}: utility designs rules and downloads rules to inverters. \emph{Right:} Inverters follow control rules fed by local and/or remote data.}
	\label{fig:diagram}  
    \end{figure*}
	
	Since the constraints in \eqref{eq:rpcfnreg} are enforced for the scenario data, the learned rules do not necessarily satisfy these constraints for all $\bz_{n,s}$ with $s\notin \{1,\ldots,S\}$. This limitation appears also in scenario-based and chance-constrained designs~\cite{Ayyagari17}. Once a control rule is learned, in real-time $t$, it can be heuristically projected within $[-\bar{q}_{n,t}^g,+\bar{q}_{n,t}^g]$ as
	\[\mcP_{\bar{q}_{n,t}^g}\left[q_{n,t}^g\right]:=\max\left\{\min\left\{q_{n,t}^g,\bar{q}_{n,t}^g\right\},-\bar{q}_{n,t}^g\right\}.\]
	

\subsection{Implementing reactive control rules}\label{subsec:implementation}
Our control scheme involves four steps; see also Fig.~\ref{fig:diagram}:
	\renewcommand{\labelenumi}{\emph{T\arabic{enumi})}}
	\begin{enumerate}
		\item The utility collects scenario data $\bz_{n,s}$ for all $n$ and $s$.
		
		\item The utility designs rules by solving \eqref{eq:rpcfnreg}; see Section~\ref{sec:svm}.
		
		\item Each inverter $n$ receives $S+1$ data $(\ba_n,b_n)$ from the utility, which describe $f_n$.
		
		\item Over the next $30$ minutes and at real time $t$, each inverter $n$ will be collecting $\bz_{n,t'}$ and applying the rule
		\begin{equation}\label{eq:rule}
		\mcP_{\bar{q}_{n,t'}^g}\left[ \sum_{s=1}^S K_n(\bz_{n,t'},\bz_{n,s}) a_{n,s} +b_n\right].
		\end{equation}
	\end{enumerate}
	
The aforesaid process is explicated next. Regarding \emph{T1)}, scenario data should be as representative as possible for the grid conditions anticipated over the following $30$-min control period. One option would be to use load and solar generation forecasts. A second option would be to use historical data from the previous day and same time, if they representative of today's conditions. A third alternative would be to use the most recent grid conditions known to the utility. For example, if smart meter data are collected every $30$~min anyway, they can be used in lieu of forecasts for the next control period.

The numerical tests of Section~\ref{sec:tests} adopt the third option and use the minute-based grid conditions observed over the last $30$-minutes as $S=30$ scenarios to train the inverter rules for the upcoming $30$-minute interval. Obviously, the number of training scenarios $S$ does not have to coincide with the length of the control period measured in minutes. These two parameters relate to loading conditions; feeder details; availability and quality of scenario data; communication and computational resources. Selecting their optimal values goes beyond the scope of this work.

During \emph{T4)}, inverter $n$ has already received $(\ba_n,b_n)$ and $\{\bz_{n,s}\}_{s=1}^S$ during \emph{T3)}. Each $\bz_n$ may consist of local data and a few active flow readings collected from major lines or transformers. If the entries of $\bz_n$ are all local, the rule can be applied with no communication. Otherwise, the non-local entries of $\bz_n$ have to be sent to inverter $n$. If non-local inputs are shared among inverters, broadcasting protocols can reduce the communication overhead.

\begin{remark}\label{re:comparison}
Suppose each inverter $n$ knows the training data $\bz_{n,s}$ for $s\in\mcS$. Function $f_n$ can be described in two ways: Either through \eqref{eq:fn} using the data described under T3); or through \eqref{eq:qfunlin}--\eqref{eq:qfunnlin} via $\bw_n$. For the second way, vector $\bw_n$ has $M_n$ entries in the linear case and $\Phi_n$ entries in the nonlinear case. For the linear case, if $M_n<S+1$, representing $f_n$ through \eqref{eq:qfunlin} by $\bw_n$ is more parsimonious. Representation \eqref{eq:fn} becomes advantageous only when $\Phi_n\gg S+1$ under the nonlinear case. 
\end{remark}

	\section{Support Vector Reactive Power Control}\label{sec:svm}	
	This section converts \eqref{eq:rpcfnreg} to a vector minimization and explores different options for $\Delta$. From \eqref{eq:rtn}, the output of inverter $n$ across all $S$ scenarios is $\bK_n\ba_n+b_n\bone$. Then, the apparent power constraints in \eqref{eq:rpcfnreg} can be written as
	\begin{equation}\label{eq:app}
	-\bar{\bq}_n^g\leq \bK_n\ba_n+b_n\bone \leq \bar{\bq}_n^g,~ \forall n
	\end{equation}
	where $\bar{\bq}_n^g := [\bar{q}_{n,1}^g ~ \cdots ~ \bar{q}_{n,S}^g]^\top$. Moreover, the vector of voltage deviations can be expressed as
	\begin{align}\label{eq:Xq+y}
	\bX\bq^g_s+\by_s&=\bX\left(\sum_{n=1}^N \be_n q_{n,s}^g\right) + \by_s\nonumber\\
	&=\sum_{n=1}^N \bx_n\be_s^\top\bK_n\ba_n+\sum_{n=1}^N b_n\bx_n +\by_s
	\end{align}
	where $\bx_n$ is the $n$-th column of $\bX$. Substituting \eqref{eq:norms} and \eqref{eq:app}--\eqref{eq:Xq+y}, the optimization in \eqref{eq:rpcfnreg} can be posed as a second-order cone program (SOCP) over $\{\ba_n\}_{n\in\mcN}$ and $\bb$. 

Nonetheless, solving \eqref{eq:rpcfnreg} with $\Delta=\Delta_\epsilon$ yields optimal $\ba_n$'s with several non-zero entries. This means that to describe rule $n$ by \eqref{eq:rule}, the utility needs to communicate the entire vector $\ba_n$ during \emph{T3)}. If scenarios $\{\bz_{n,t}\}_{t=1}^T$ are not known by the inverter, they have to be communicated along with $\ba_n$ as well. The number of scenarios $T$ may be large when learning rules under complex feeder setups. A related approach for minimizing a convex combination of $\Delta_s$ and power losses has been suggested in the conference precursor of this work~\cite{G-SIP18}, but inherits the same difficulty of non-sparse $a_{n,t}$'s.

Inspired by support vector machines (SVM), we engineer $\Delta$ to obtain inverter rules described by possibly fewer scenarios: Promoting sparse $\ba_n$'s alleviates the communication overhead during step \emph{T3)}. To this end, we put forth the cost
	\begin{equation}\label{eq:dvt}
	\Delta_\tau(\bq^g;\by):=\left[\|\bX\bq^g+\by\|_2\right]_\tau
	\end{equation}
	for some $\tau>0$. If scenario $s$ yields a vector of voltage deviations $\bX\bq^g_s+\by_s$ with $\ell_2$-norm smaller than $\tau$, this scenario incurs no cost. If $\|\bX\bq^g_s+\by_s\|_2>\tau$, the voltage regulation penalty grows with $\|\bX\bq_s^g+\by_s\|_2$. The cost in \eqref{eq:dvt} can be expressed as an SOCP over the slack variable $d$
	\begin{align*}
	\Delta_\tau(\bq^g;\by):=\min_{d\geq 0}\left\{d:\|\bX\bq^g+\by\|_2\leq d+\tau\right\}.
	\end{align*}
	Applying the same epigraph trick for the function norms, problem \eqref{eq:rpcfnreg} can be solved as the SOCP
	\begin{subequations}\label{eq:kvr}
		\begin{align}
		\min~&~\frac{1}{S}\bd^\top\bone + \mu \bgamma^\top\bone\label{eq:kvr:cost}\\
		\mathrm{over}~&~\{\bq_s^g\},\{\ba_n\}, \bb, \bd\geq \bzero, \bgamma \label{eq:kvr:vars}\\
		\mathrm{s.to}~&~\eqref{eq:app},\eqref{eq:Xq+y}\label{eq:kvr:app}\\
		~&~\|\bK_n^{1/2}\ba_n\|_2\leq \gamma_n,\quad \forall n \label{eq:kvr:norms}\\
		~&~ \|\bX\bq_s^g+\by_s\|_2\leq d_s+\tau,\quad \forall s \label{eq:kvr:svm}
		\end{align}
	\end{subequations}
where $\bd:=[d_1~\cdots~d_S]^\top$ and $\bgamma := [\gamma_1~\cdots~\gamma_N]^\top$. The variables $\bq_s^g$ can be eliminated using the substitutions of \eqref{eq:Xq+y}. Solving \eqref{eq:kvr} takes $\mcO\left(N^{3.5}T^3\right)$ operations with interior point-based solvers~\cite{SOCP_Lobo}. However, the advantage of inverter control rules is that \eqref{eq:kvr} is not solved in real time. If standard interior point-based solvers are not scalable to larger grids, one may resort to (distributed) first-order algorithms; warm-start initializations; and cutting-plane methods.

The coefficients $\ba_n$'s minimizing \eqref{eq:kvr} enjoy two types of sparsity, across inverters and across scenarios. To explain the first type of sparsity, express the second summand in the cost of \eqref{eq:kvr} as $\mu\bgamma^\top\bone = \mu\sum_{n=1}^N \|\bK_n^{1/2}\ba_n\|_2$. Having these non-squared $\ell_2$-norms in the objective promotes block sparsity across $n$, in the sense that for larger $\mu$, some vectors $\ba_n$ may be set to zero. This effect is a direct consequence of block-sparse solutions encountered in group Lasso (G-Lasso)-formulations; see~\cite{Ravikumar09}, \cite{BaGia13}, \cite{KZG14}. All inverters receive a reactive power setpoint $b_n$, but if the optimal $\ba_n$ becomes zero, inverter $n$ will not be changing its reactive injection in real-time. One may drop the intercept $b_n$ from the control rule of \eqref{eq:qfun} and the optimization of \eqref{eq:kvr}, and modify the feature vector as 
\begin{equation}\label{eq:augment}
\bz_n'=[1~\bz_n^\top]^\top. 
\end{equation}
Thus, obtaining $\ba_n=\bzero$ from \eqref{eq:kvr} enables inverter selection. 

The next proposition studies the second type of sparsity; see the appendix for a proof.

\begin{proposition}\label{pro:sparsity}
Consider \eqref{eq:rpcfnreg} with $\Delta=\Delta_\tau$ and its minimizer in \eqref{eq:fn}. If $\|\bX \bq_s^g +\by_s\|_2< \tau$ for scenario $s$ at the optimum, then $a_{n,s}=0$ for every inverter $n$ with $|q_{n,s}^g|<\overline{q}_{n,s}^g$. 
\end{proposition}

Proposition~\ref{pro:sparsity} explains how $\Delta_\tau$ promotes block sparsity across $s$: If scenario $s$ does not experience severe voltage violations, the corresponding coefficients $a_{n,s}$ will be zero for all inverters $n$ that have not reached their apparent power limit. Block sparsity across time identifies non-critical scenarios. Phrased in the SVM context, the so-termed `support vectors' here correspond to scenarios with significant voltage deviations. Larger values of $\tau$ effect fewer critical scenarios. 

These two forms of sparsity offer communication savings since the related $(a_{n,s},\bz_{n,s})$ do not need to be communicated to inverters. This enables training the rules for larger number of scenarios $S$ at the same communication overhead. Note that for fixed $(\mu,\tau)$, the sparsity of $\ba_n$'s depends on the training data $\by_s$'s as well. If a particular sparsity goal is to be met, the utility has to solve \eqref{eq:kvr} repeatedly for various values of $\mu$ and $\tau$. Such computations can be significantly sped up by initializing an optimization algorithm for one value of $\tau$ to the minimizer obtained using the previous value of $\tau$~\cite[Sec.~18.4]{Hastie}; however, such techniques will not be pursued here.
	
Different from $\Delta_\tau$, cost $\Delta_\epsilon$ is not expected to yield as sparse $\ba_n$'s. The next claim (proved in the appendix) explains that even if a single bus experiences voltage deviation larger than $\epsilon$ for scenario $s$, then $a_{n,s}\neq 0$ for all $n$. In other words, a voltage violation at a single bus for scenario $s$ renders this scenario critical for \emph{all} inverter rules.

\begin{proposition}\label{pro:sparsity2}
Consider \eqref{eq:rpcfnreg} with $\Delta=\Delta_\epsilon$ and its minimizer in \eqref{eq:fn}. If $\|\bX \bq_s^g +\by_s\|_{\infty}> \epsilon$ for scenario $s$ at the optimum, then $a_{n,s}\neq 0$ for all $n$. 
\end{proposition}

	\section{Numerical Tests}\label{sec:tests}
	The novel inverter rules were tested on the IEEE 123-bus feeder~\cite{testfeeder}, converted to a single-phase grid as described in~\cite{GLTL12}. Residential load and solar data were extracted from the Pecan Street dataset as delineated next~\cite{pecandata}. Minute-sampled active load and solar generation data were collected for June 1, 2013 between 8:00--16:00. We downloaded data from the first 123 Pecan Street nodes, after excluding nodes with empty data records. Regarding solar generation, unless stated otherwise, $75\%$ of the buses had solar generation by excluding nodes with bus indexes that are multiples of $4$. 
	
	
	 Load data were scaled on a per bus basis so that their daily peak values matched $150\%$ of the benchmark load. Since the Pecan Street data  included only active power, we drew lagging power factors uniformly at random within $[0.9,0.95]$ for each bus and kept them fixed across time. The scaling factors for active loads were also used for scaling solar data. To allow for reactive power compensation even at peak solar irradiance, inverters were over-sized by $10\%$ providing an apparent power capacity of $\bar{s}_n^g=1.1\bar{p}_n^g$ for all $n$; see~\cite{Turitsyn11}. 
	
	Our numerical tests included six control schemes:\\
	\hspace*{1em}\emph{C1)} The optimal reactive injections computed by \eqref{eq:vr} on a per-minute basis;\\
	\hspace*{1em}\emph{C2)} The optimal reactive injections computed by \eqref{eq:vr} on a per-minute basis assuming a $2$-minute communication delay;\\
	\hspace*{1em}\emph{C3)} The fixed Watt-VAR control rules of \cite[(12)--(14)]{Turitsyn11};\\
	\hspace*{1em}\emph{C4)} The rules of \eqref{eq:rpcfnreg} for linear kernels and $\Delta=\Delta_\tau$;\\
	\hspace*{1em}\emph{C5)} The rules of \eqref{eq:rpcfnreg} for Gaussian kernels and $\Delta=\Delta_\tau$;\\
	\hspace*{1em}\emph{C6)} The rules of \eqref{eq:rpcfnreg} for linear kernels and $\Delta=\Delta_\epsilon$; and\\
	\hspace*{1em}\emph{C7)} The rules of \eqref{eq:rpcfnreg} for Gaussian kernels and $\Delta=\Delta_\epsilon$.
	
The input $\bz_n$ to inverter $n$ consisted of local data as in \eqref{eq:input}. Each entry of $\bz_n$ was centered by its daily mean and normalized by its daily standard deviation. To avoid rank deficiency, we added $10^{-3}\cdot\bI_S$ to all kernel matrices. 

 Schemes \emph{C1),C2)} were solved using SDPT3 and YALMIP with MATLAB~\cite{yalmip,SDPT3}. Schemes \emph{C4)--C7)} were solved by invoking the MOSEK solver directly through MATLAB~\cite{mosek}. Tests were run on a 2.4 GHz Intel Core i5 laptop with 8 GB RAM. The average running time for solving~\eqref{eq:rpcfnreg} with $T=30$ is given in Table~\ref{table:time}. It should be emphasized that although the control rules were designed using the LDF grid model, the voltage deviations experienced by all control rules were tested using the full AC model. 

	\begin{table}[t]
		\begin{center}
			\caption{Running Time for Solving~\eqref{eq:rpcfnreg} with $T=30$} 
			\label{table:time}
			\begin{tabular}{|c|c|c|c|c|}\hline
				\hline
				& \emph{C4)} & \emph{C5)} & \emph{C6)} & \emph{C7)} \\ 
				\hline\hline
				Running time [min] & $0.21$ & $0.45$ & $0.96$ & $1.99$\\ 
				\hline \hline
			\end{tabular}
		\end{center}
	\end{table}
		
During training, we used $T=30$ scenarios to learn the SVM-based control rules of  \emph{C4)--C7)}. These scenarios comprised the load and solar data observed during the last $30$~minutes. During validation, the inverter control rules were tested over the following $30$ minutes. Parameters $\mu$ and $\gamma$ were selected via $5$-fold cross-validation. The ranges of $\tau$ and $\epsilon$ were empirically chosen to yield an average communication overhead similar to the one needed by the affine rule of \eqref{eq:qfunlin} as discussed under Remark~\ref{re:comparison}: An affine rule is described by $M_n+1=4$ data per inverter. If only $10\%$ of the entries of $\ba_n$ are nonzero, then communicating $(\ba_n,b_n)$ entails sending $0.1\cdot S+1=0.1\cdot 30+1=4$ data as well. The sparsity of $\ba_n$'s depends on input data along with the values of $(\tau,\mu)$ or $(\epsilon,\mu)$. These parameters were set so that $\ba_n$'s had  $10\%$ nonzero entries on the average across time and buses.

\begin{figure}[t]
		\centering
		\includegraphics[scale=0.24]{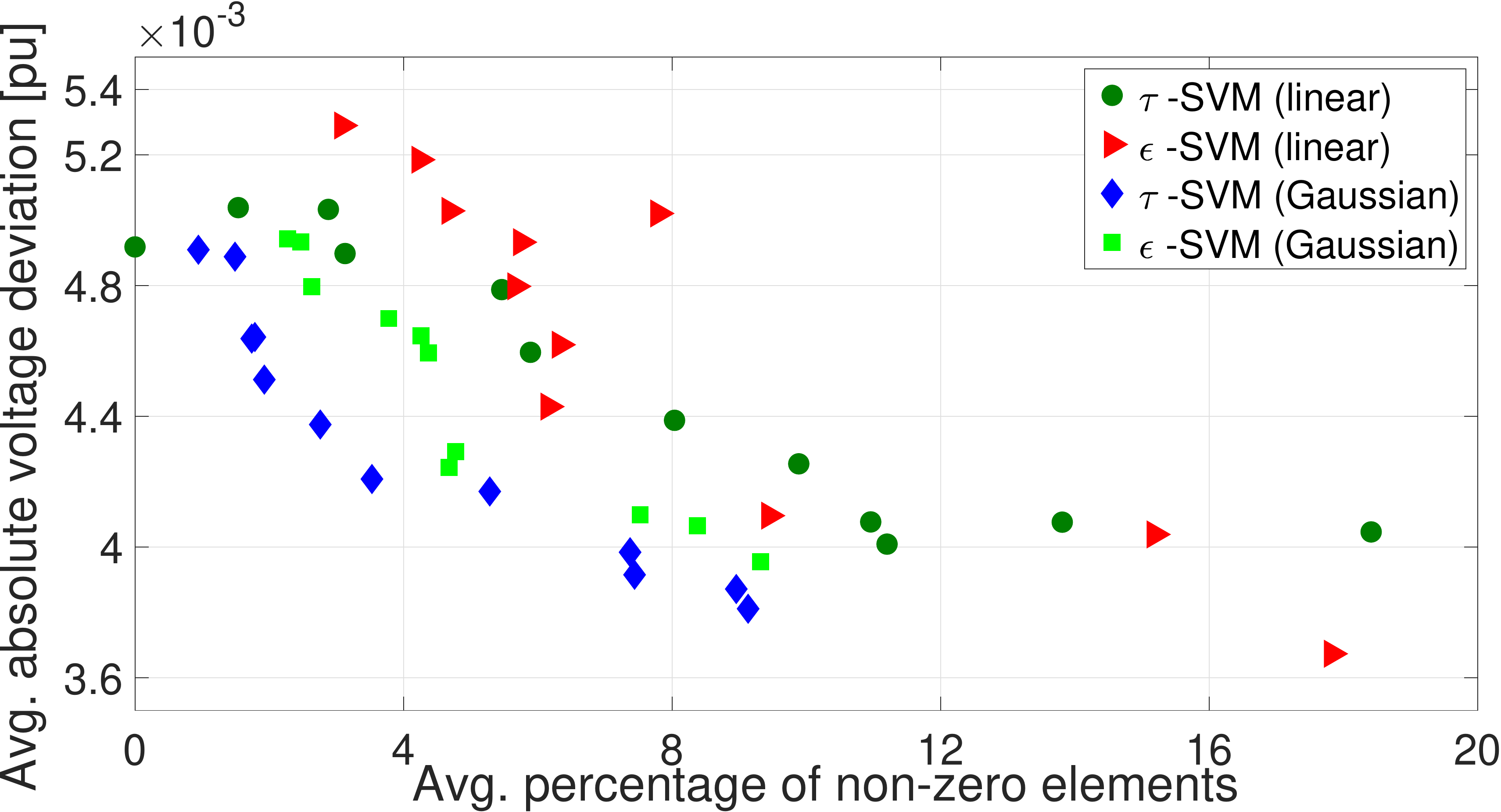}
		\caption{Average of absolute voltage deviation vs. sparsity for \emph{C3)--C6)}.}
		\label{fig:avg-performance}  
	\end{figure}

We next explored the trade-off between voltage deviation and the sparsity of $\ba_n$'s for \emph{C4)--C7)}. The expectations from this test were two: \emph{i)} voltage deviations are expected to increase for sparser $\ba_n$'s; \emph{ii)} schemes \emph{C4)} and \emph{C5)} should exhibit improved sparsity over \emph{C6)} and \emph{C7)}. To validate these hypotheses, we recorded the voltage deviations for $10$ values of $\tau$ and $\epsilon$ for \emph{C4)--C7)}. The average absolute voltage deviation and the average percentage of non-zero coefficients were calculated over the day and across buses, and are shown in Figure~\ref{fig:avg-performance}. From Figure~\ref{fig:avg-performance}, the value of $\tau$ yielding a sparsity of roughly $11\%$ is $\tau=0.001$. Figure~\ref{fig:avg-performance} reveals three important points. First, voltage deviations increase as $\ba_n$'s become sparser as expected. Second, for a given sparsity in $\ba_n$'s, the rules obtained by $\Delta_\tau$ exhibit smaller voltage deviations compared to the rules obtained by $\Delta_{\epsilon}$. Because of this, we focus on the performance of \emph{C4)--C5)} for the rest of this section. Third, the Gaussian kernel-based rules attained lower voltage deviations than the related linear kernel-based rules.
	
\begin{figure}[t]
	\centering
	\includegraphics[scale=0.24]{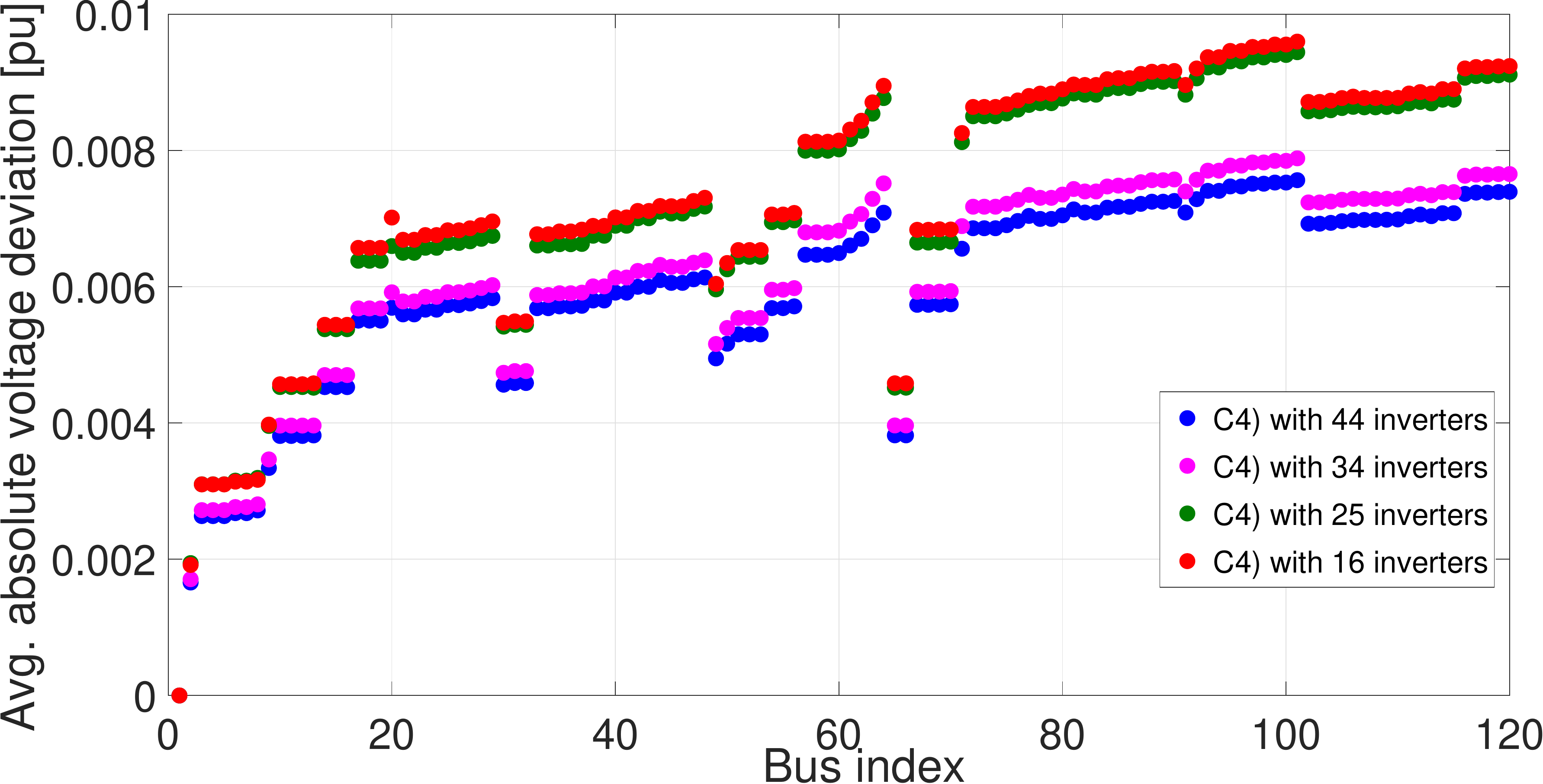}
	\caption{Maximum absolute voltage deviation over time for $75\%$ penetration, obtained by the Gaussian SVM-based rules trained for $\Delta=\Delta_\tau$.}
	\label{fig:inverters}  
    \end{figure}

We next tested the effect of $\mu$ on inverter selection and voltages. Larger values of $\mu$ are expected to set more $\ba_n$'s to zero. To eliminate the inverters with $\ba_n=\bzero$, the parameter $b_n$ was appended in $\ba_n$ as delineated in \eqref{eq:augment}. For a fixed value of $\tau=0.001$, for scheme~\emph{C4)}, the values of $\mu$ were obtained using cross-validation across the day. The control rules were designed again using $4$ different values of $\mu$. As expected, by increasing the value of $\mu$, the number of all-zero $\ba_n$'s and the corresponding voltage deviations were increased. Figure~\ref{fig:inverters} depicts the absolute voltage deviation averaged over time for each inverter. Notice that the values of $\tau$ and $\mu$ were kept fixed, although the training data $\by_s$'s varied across the day. Due to this, the reported sparsity in Figure~\ref{fig:avg-performance} is the average sparsity across time and inverters. Moreover, the number of inverters in Figure~\ref{fig:inverters} is the average number of activated inverters across the day. Even though the values of $\mu$ and $\tau$ can be adjusted on a $30$-min basis to meet specific sparsity requirements, we chose to keep them fixed to simplify the exposition. In fact, the rest of this section reports the worst-case instead of average voltage deviations across time and for each bus.

\begin{figure}[t]
		\centering
		\includegraphics[scale=0.24]{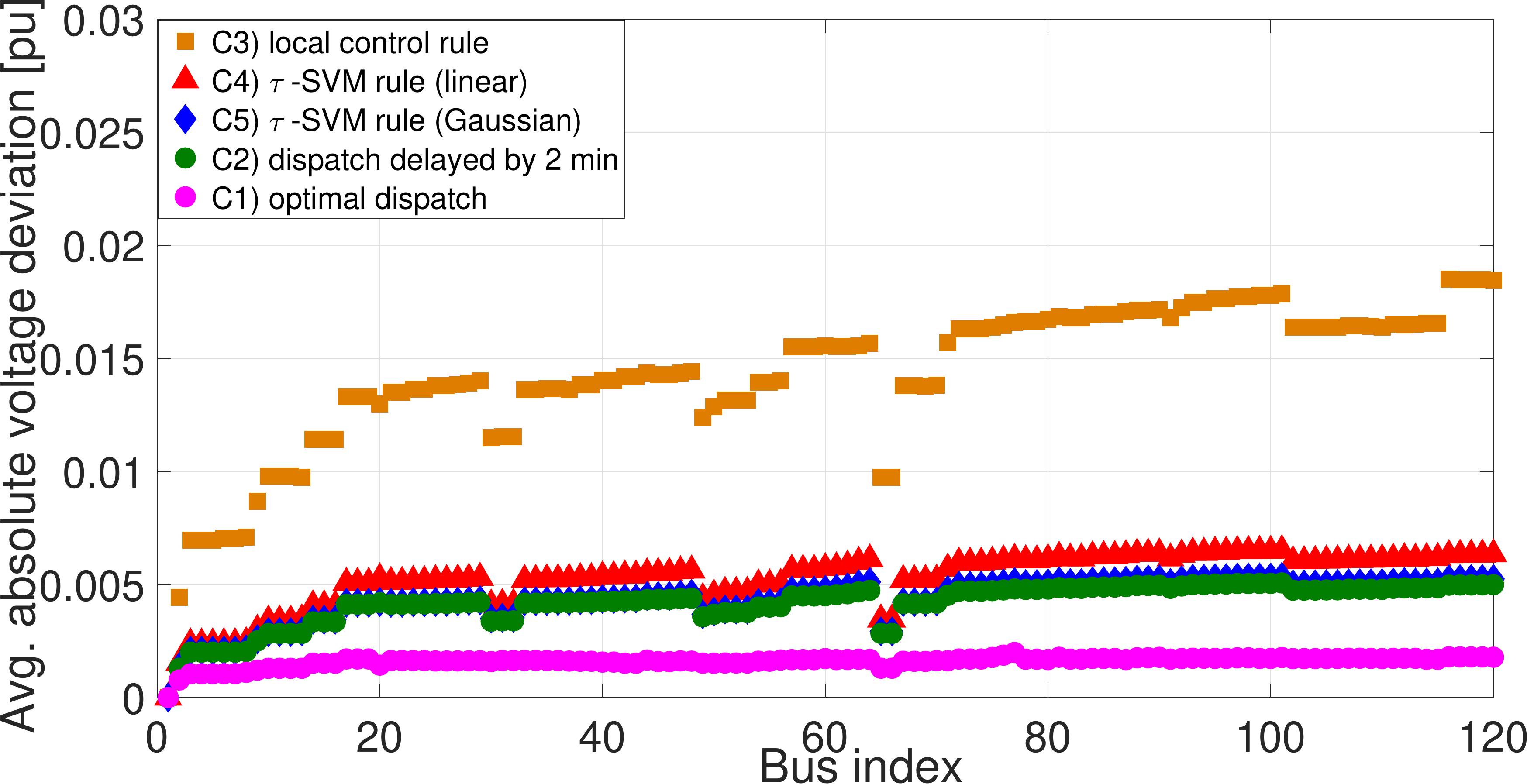}
		\caption{Absolute voltage deviation averaged over time for $75\%$ penetration, obtained by the SVM-based rules trained for $\Delta=\Delta_\tau$.}
		\label{fig:dv70}  
	\end{figure}

\begin{figure}[t]
	\centering
	\includegraphics[scale=0.24]{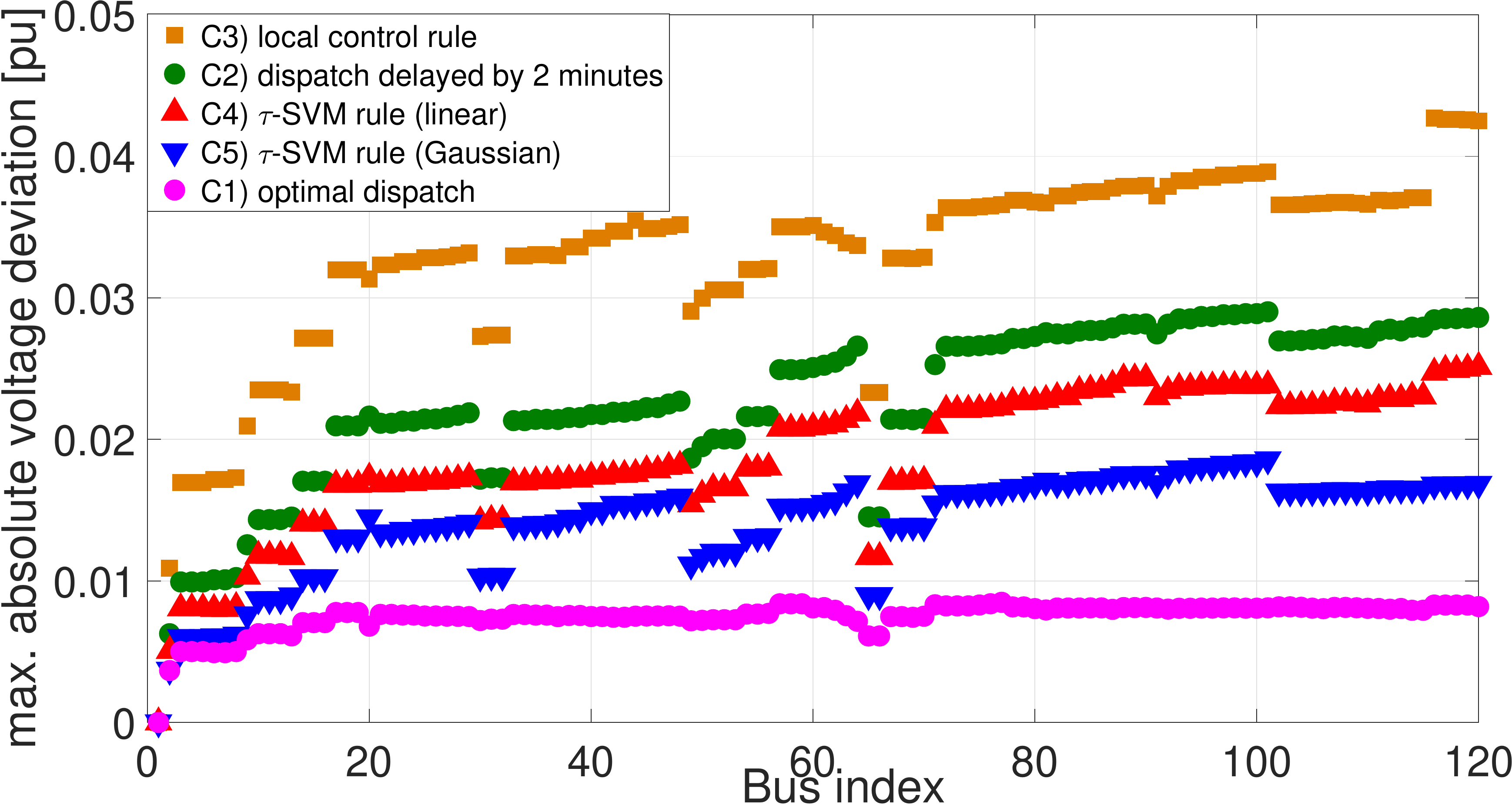}
	\caption{Maximum of absolute voltage deviation  over time for $75\%$ penetration, obtained by the SVM-based rules trained for $\Delta=\Delta_\tau$.}
	\label{fig:max70}  
\end{figure}

We next compared the proposed SVM-based control rules against the alternative schemes of \emph{C1)--C3)}. To this end, voltage deviations were calculated between 8:00--16:00 for schemes \emph{C1)--C5)}. Figures~\ref{fig:dv70} and~\ref{fig:max70} demonstrate the average and the maximum voltage deviations over the test period. It can be observed from both figures that the Gaussian SVM-based rule performs better than \emph{C1)--C3)} due to its ability to capture non-linear behaviors. Although \emph{C3)} needs no communication, it violates the ANSI-C.84.1 standard voltage constraints. Furthermore, despite the high communication needed, scheme \emph{C2)} shows no superiority in performance over \emph{C5)} and corroborates the need for real-time response to system inputs. 

\begin{figure}[t]
		\centering
		\includegraphics[scale=0.24]{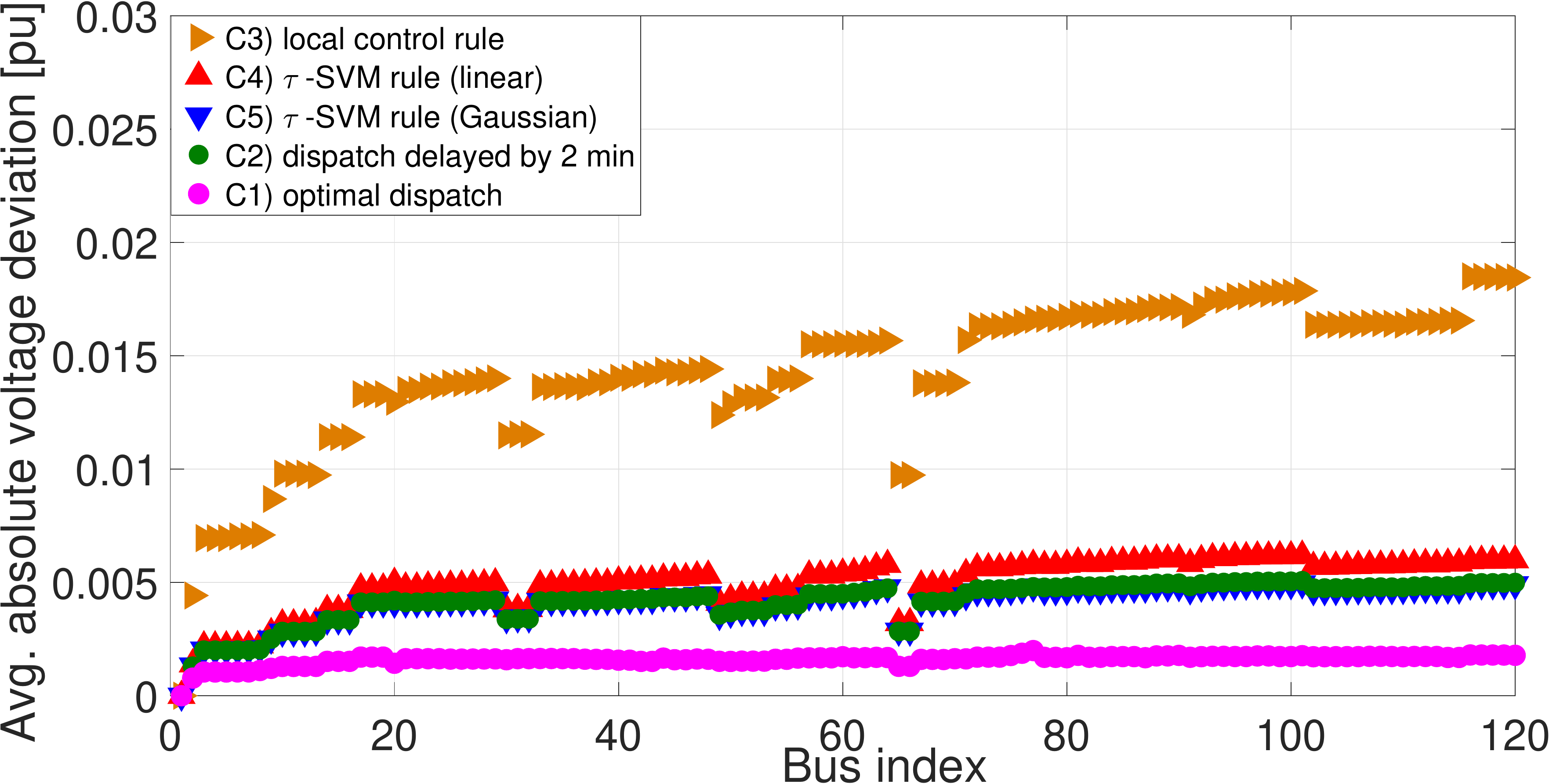}
		\caption{Absolute voltage deviation averaged over time for $75\%$ penetration with remote inputs, and the SVM-based rules trained for $\Delta_\tau$.}
		\label{fig:dv70_global}  
	\end{figure}

\begin{figure}[t]
	\centering
	\includegraphics[scale=0.24]{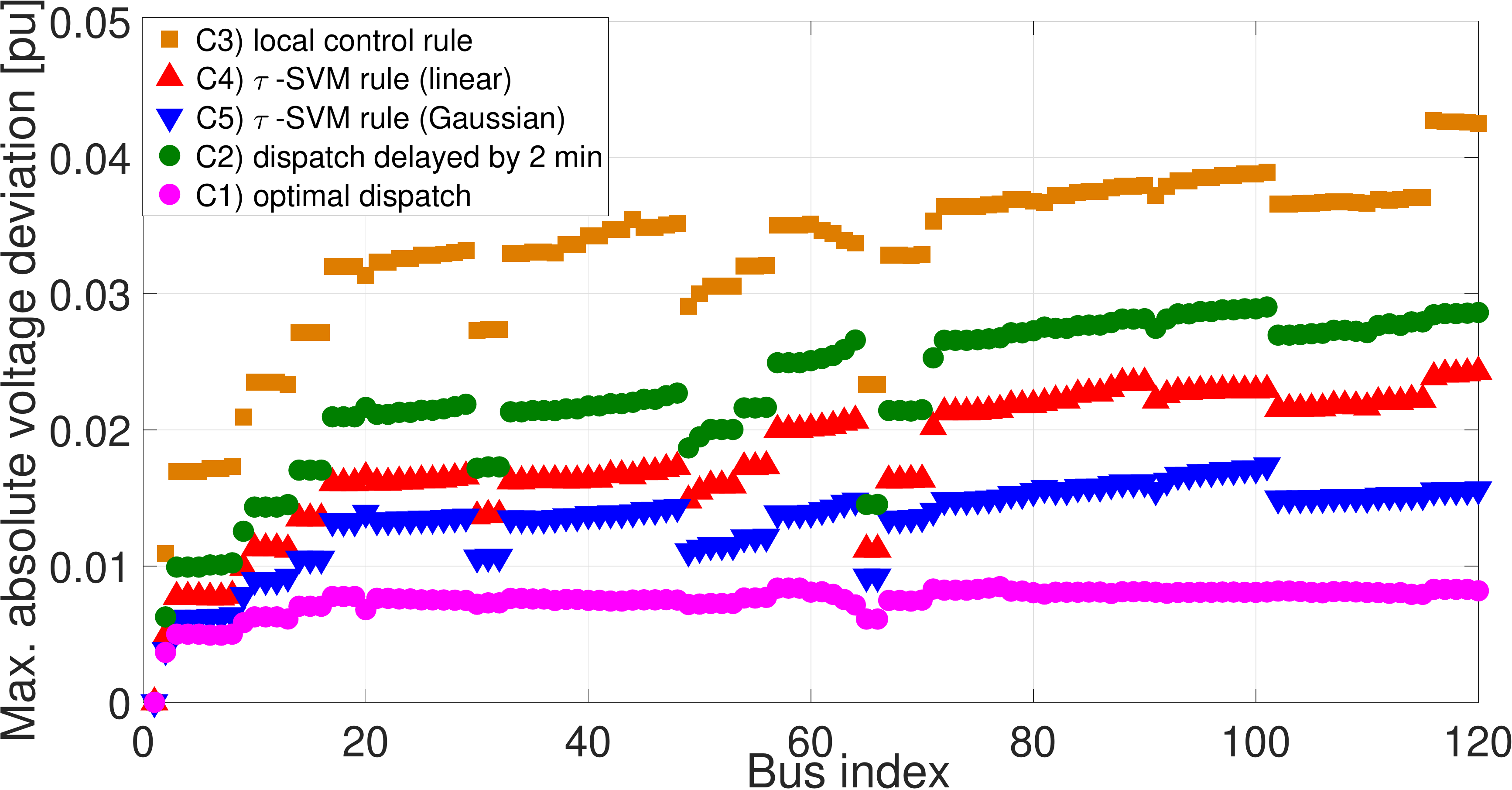}
	\caption{Maximum absolute voltage deviation over time for $75\%$ penetration obtained with remote inputs and the SVM-based rules trained for $\Delta=\Delta_\tau$.}	
	\label{fig:max75_global}  
\end{figure}

In all previous tests, the rules were fed with locally recorded data. To evaluate the advantage of adding remote control inputs, we appended the values of active power flows on the lines feeding buses $1$, $16$, and $51$, to all input vectors $\bz_n$. The daily maximum and the average voltage deviations attained by \emph{C1)--C5)} are depicted in Figures~\ref{fig:dv70_global} and~\ref{fig:max75_global}, respectively. As expected, the results suggest that adding remote inputs to the rules improves the grid voltage profile at the expense of increased inter-network communication.
	
\begin{figure}[t]
	\centering
	\includegraphics[scale=0.24]{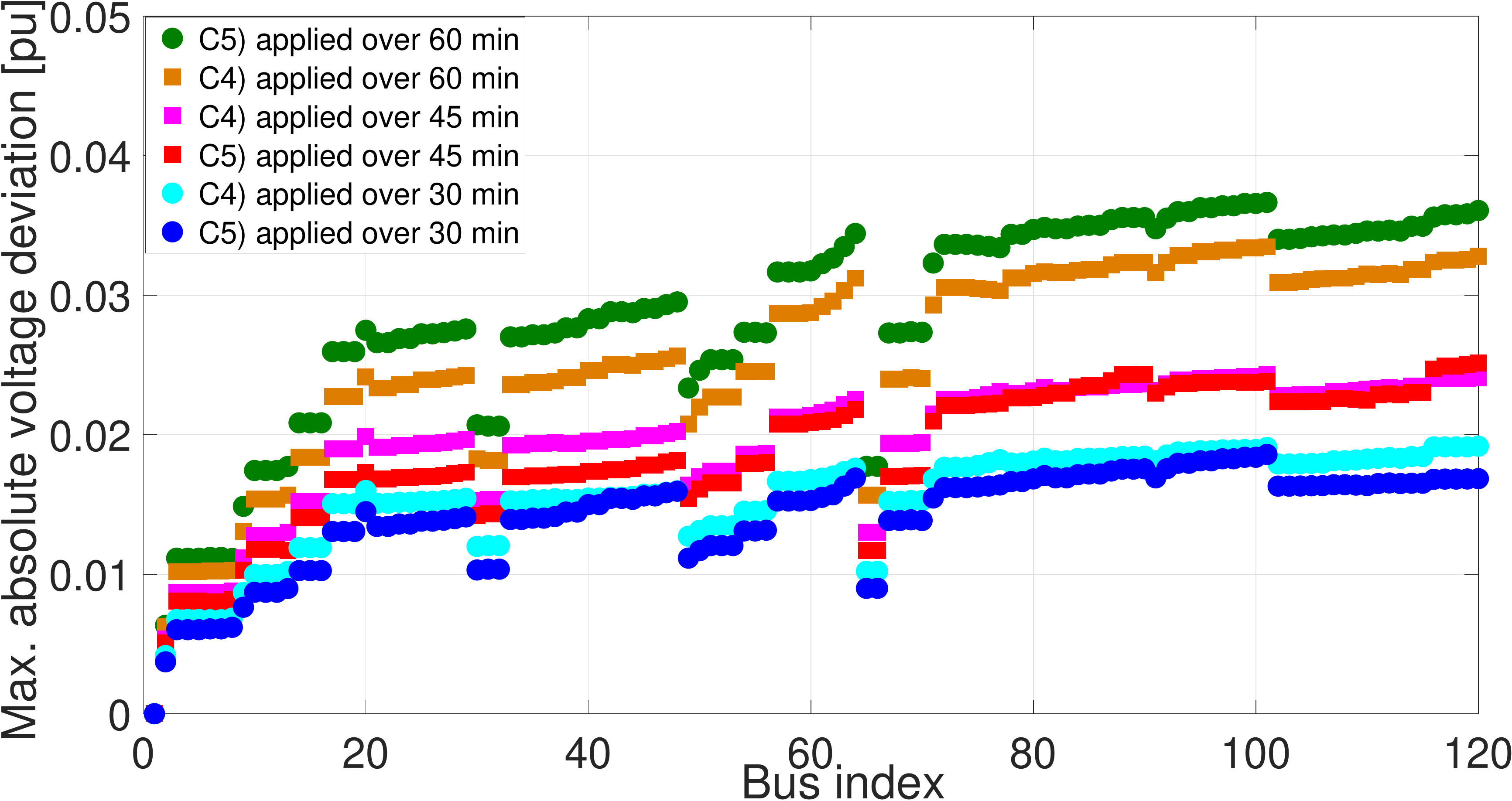}
	\caption{Maximum absolute voltage deviation over time for $75\%$ penetration obtained using the SVM-based rules applied over $30$, $45$, and $60$ minutes.}
	\label{fig:60-30}  
    \end{figure}

As mentioned in Section~\ref{subsec:implementation}, the length of the control period (in minutes) over which rules remain constant does not have to agree with the number of scenarios $S$ used for training the rules. To evaluate how the control rules perform for longer control periods, Figure~\ref{fig:60-30} compares the voltage deviations obtained by training rules using $S=30$ scenarios, but keeping them unaltered over $30$, $45$, and $60$ minutes. As expected, voltage regulation deteriorates as rules remain unchanged for longer periods.
	
\begin{figure}[t]
	\centering
	\includegraphics[scale=0.24]{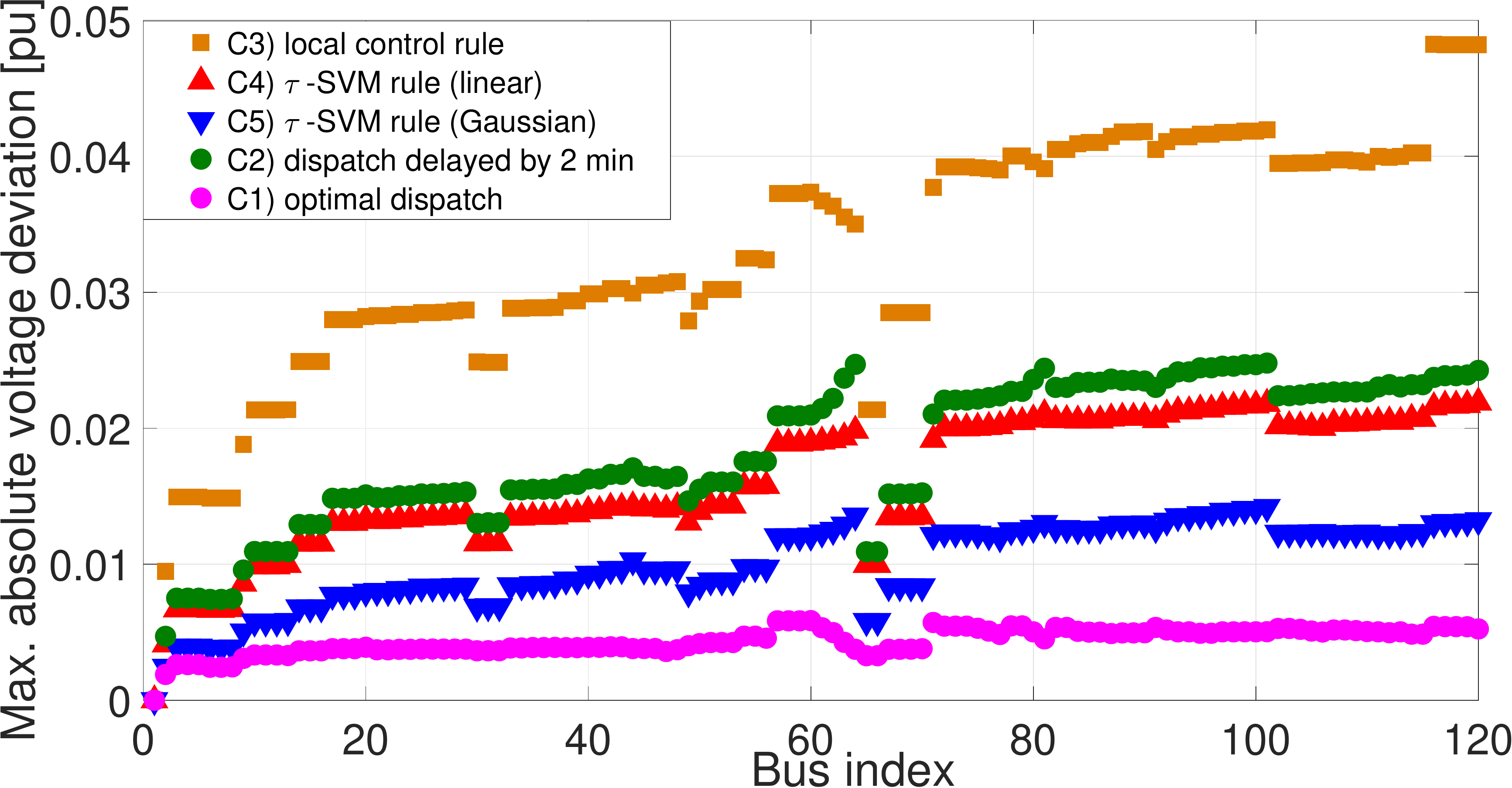}
	\caption{Maximum absolute voltage deviation over time for $50\%$ penetration obtained using the SVM-based rules trained for $\Delta=\Delta_\tau$.}
	\label{fig:max50}  
    \end{figure}

\begin{figure}[t]
	\centering
	\includegraphics[scale=0.24]{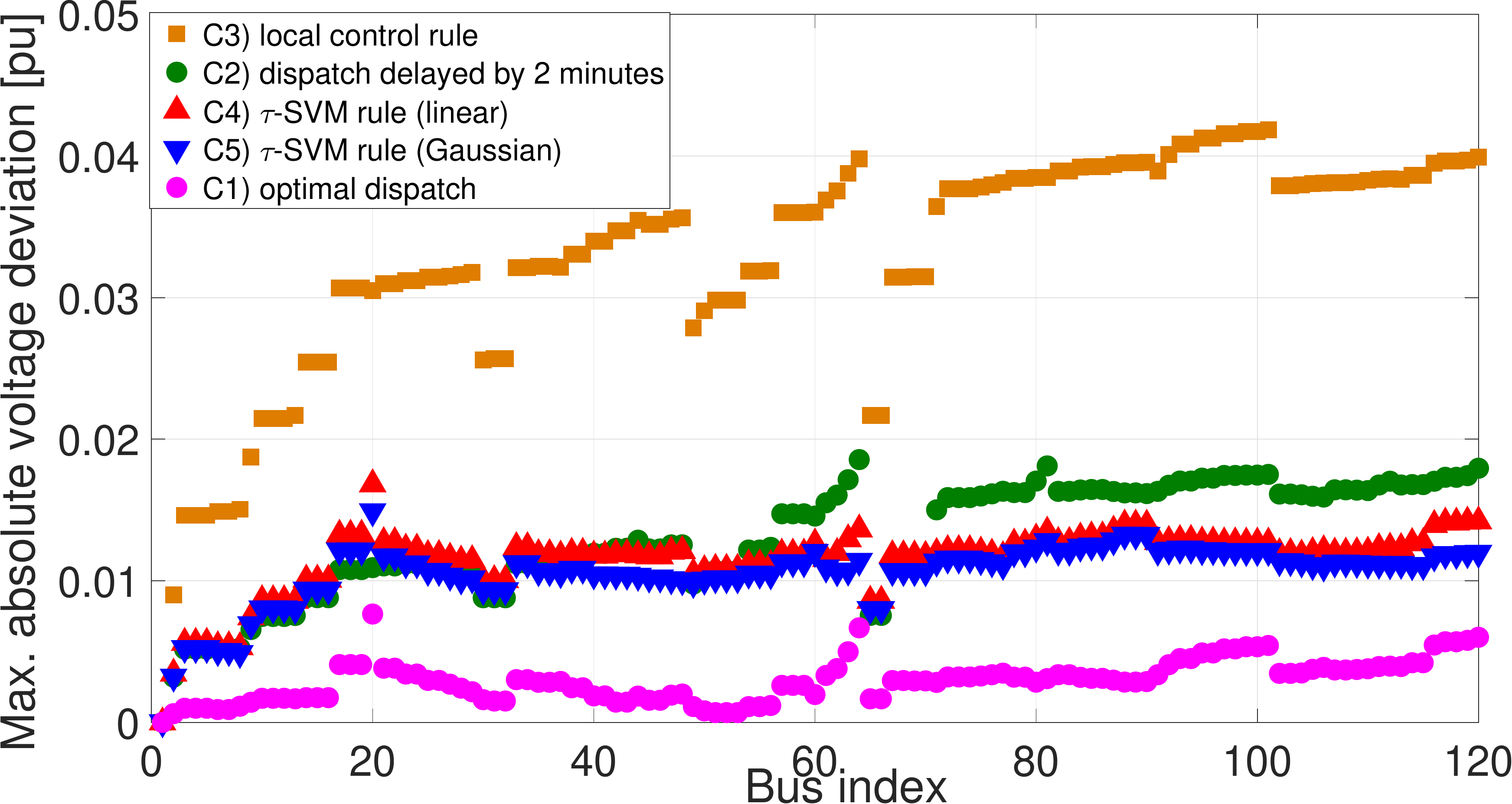}
		\caption{Maximum absolute voltage deviation over time for $25\%$ penetration obtained using the SVM-based rules trained for $\Delta=\Delta_\tau$.}
	\label{fig:max25}  
  \end{figure}

All previous tests assumed solar penetration of $75\%$. We also tested the performance of \emph{C1)--C5)} under penetrations of $50\%$ and $25\%$. To simulate $50\%$ penetration, solar generation and smart inverters were installed only in buses with even indexes. Likewise, to simulate $25\%$ penetration, we considered buses whose indexes were multiples of $4$. Figures~\ref{fig:max50} and~\ref{fig:max25} depict the attained maximum absolute voltage deviations, which apparently decrease with decreasing solar penetration. For lower penetrations, the Gaussian-based rule preserves its superior voltage profile over the other schemes.

\begin{figure}[t]
	\centering
	\includegraphics[scale=0.24]{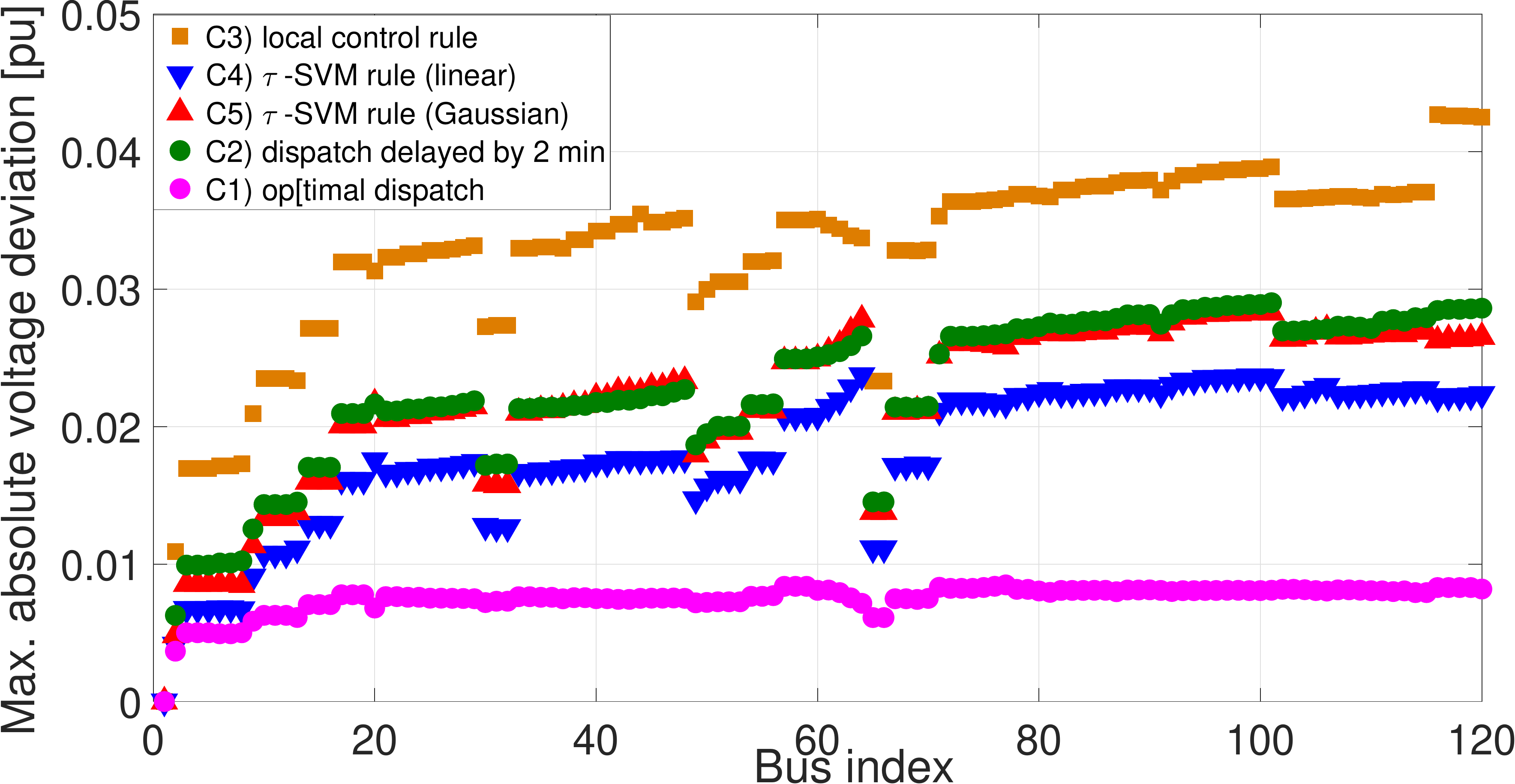}
	\caption{Max. absolute voltage deviation over time for $75\%$ penetration obtained using SVM-based rules trained for $\Delta=\Delta_\tau$ with $\tau=0.03$.}
	\label{fig:max3-75-03}  
\end{figure}

Schemes \emph{C4)} and \emph{C5)} were also tested under less communication by scaling down the sparsity in $\ba_n$'s by a factor of $10$: Voltage deviations were evaluated for $\tau=0.03$ corresponding to $1.4\%$ non-zero entries for $\ba_n$'s on the average. Figure~\ref{fig:max3-75-03} demonstrates the maximum absolute voltage deviation for \emph{C1)--C5)}. Even with fewer coefficients communicated, the voltage constraints of ANSI-C.84.1 were still satisfied. 

\begin{figure}[t]
	\centering
	\includegraphics[scale=0.24]{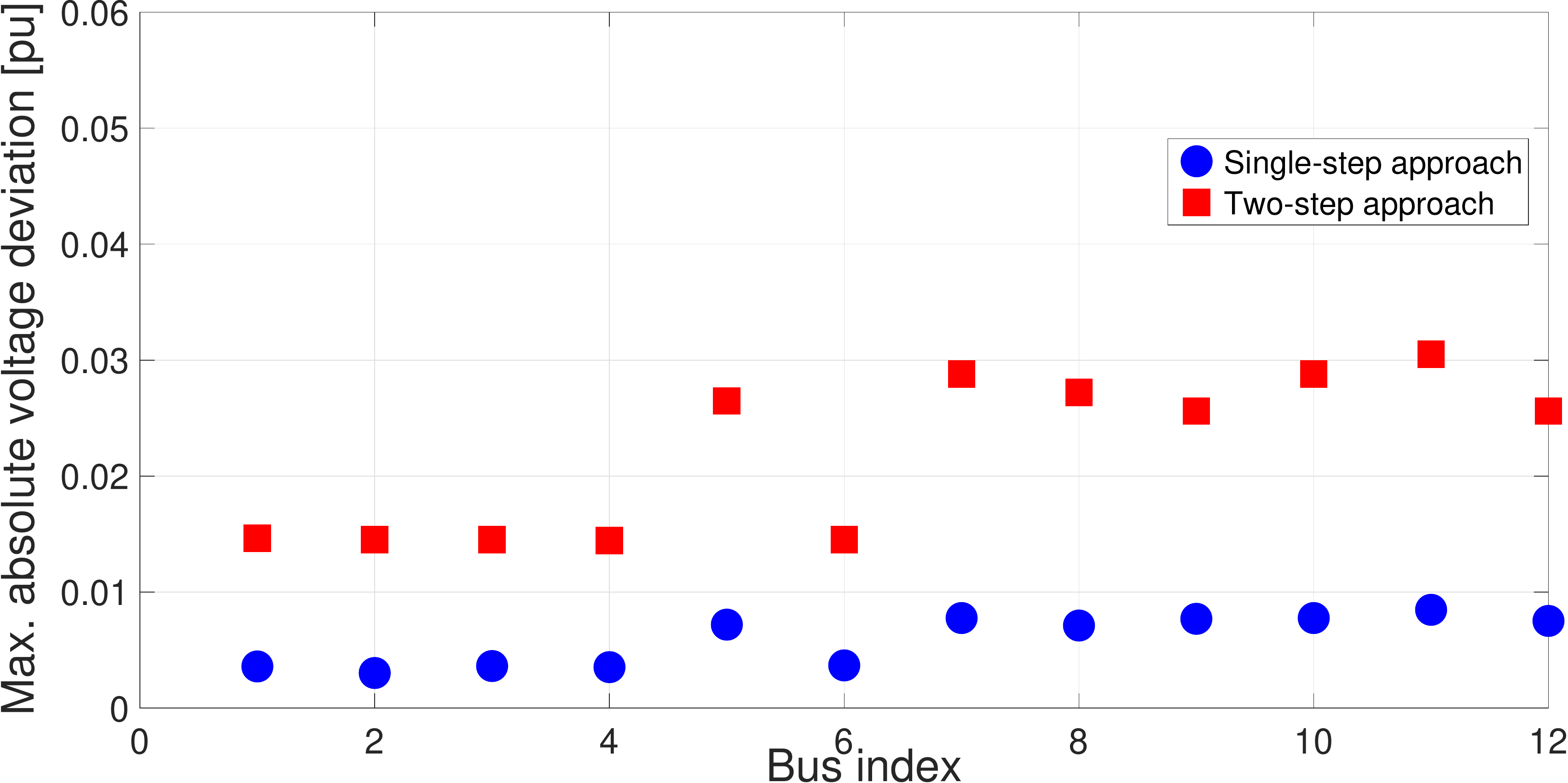}\\
	\includegraphics[scale=0.24]{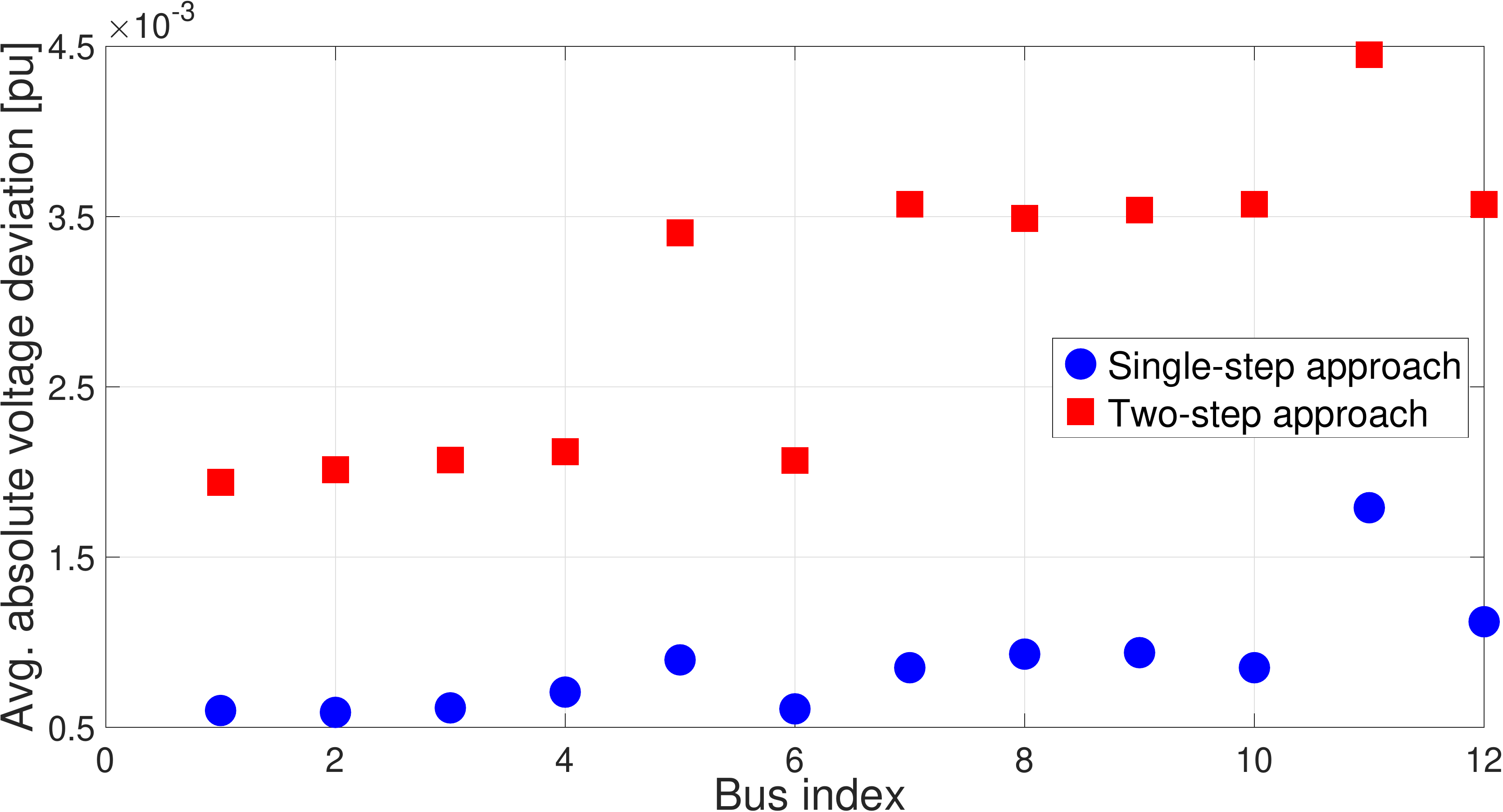}\\
	\includegraphics[scale=0.24]{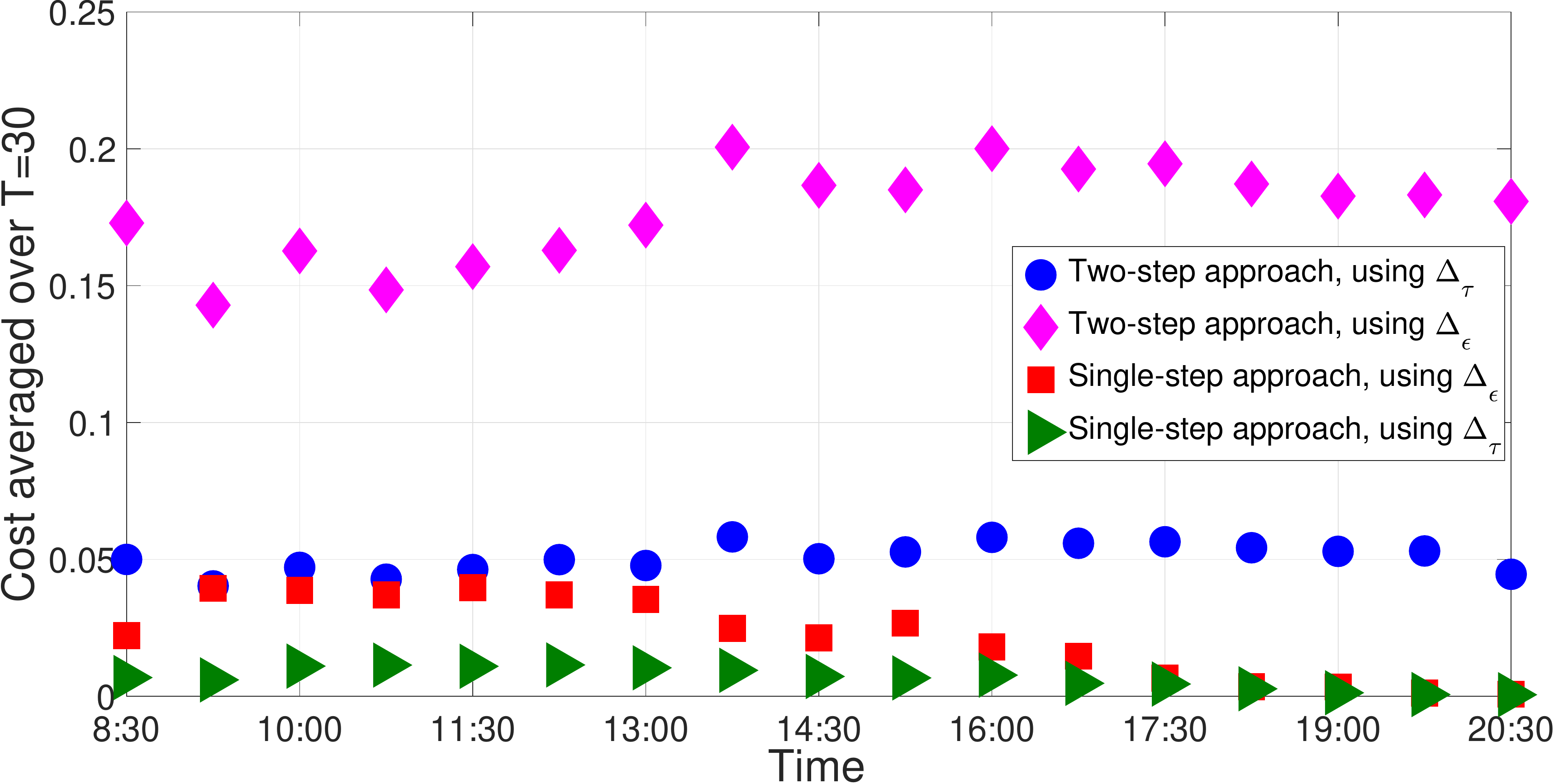} 
	\caption{Comparison between the proposed single-step learning approach (rules $R_1$), and the two-step learning approach of \cite{Dobbe19}--\cite{Kara18} (rules $R_2$).}
	\label{fig:R1vsR2}
\end{figure}

The last set of numerical tests compares the developed single-step approach with the two-step approach of \cite{Dobbe19}--\cite{Kara18}; see also Remark~\ref{re:others}. We used both approaches to design local linear control rules for the IEEE 13-bus feeder~\cite{Kersting}, under the voltage deviation cost $\Delta=\Delta_{\epsilon}$ with $\epsilon=0.001$. The top and center panels of Fig.~\ref{fig:R1vsR2} show respectively the maximum and average voltage deviation per bus computed across time. The bottom panel shows the voltage deviation cost $\Delta_\epsilon$, time-averaged per control period. The bottom panel also shows the voltage deviation cost $\Delta_\tau$ with $\tau=0.01$ attained upon training both rules using $\Delta_\tau$ instead of $\Delta_\epsilon$. Similar results were obtained for other values of $\epsilon$ and $\tau$. According to these tests, the single-step approach achieved: \emph{1)} lower maximum per-bus voltage deviations; \emph{2)} lower average per-bus voltage deviations; and \emph{3)} smaller voltage deviation costs during the operational phase.

\section{Conclusions}\label{sec:conclusions}
A novel approach for designing inverter control rules has been put forth. It relies on both data-based learning and physical grid modeling. Inverter rules are not learned independently using input/output pairs of the OPF problem. Instead, they are learned jointly by posing the related OPF problem as a multi-function learning task. Because of the way voltage deviations couple inverter outputs, the conventional support vector machine approach fails to yield sparse rule descriptions. We have engineered a voltage deviation cost to identify `support scenarios,' that is a few scenarios with non-zero coefficients for most of inverter rules. The devised control rules were tested using on a benchmark feeder using the exact AC model. The novel scheme attained superior voltage regulation performance compared to preset local rules, and oftentimes comparable performance to an optimal inverter dispatch delayed by $2$ minutes. The numerical tests have further corroborated the benefits of nonlinear rules with non-local inputs, and explored the trade-off between voltage regulation performance and sparsity. Finally, this work motivates several questions. On the implementation side, testing the novel formulations on multiphase grids along with capacitor banks, voltage regulators, and ZIP loads, is of practical interest. On the analytical side, chance-constrained formulations; studying the stability of nonlinear rules with voltages as inputs; using kernels to learn functions with constraints; and selecting non-local control inputs; are some open and interesting questions.
	
	\appendix\label{sec:appendix}
	\begin{IEEEproof}[Proof of Proposition~\ref{pro:sparsity}] 
		Consider first the linear rules of \eqref{eq:qfunlin}, for which $q_n^g(\bz_n)=\bz_n^\top\bw_n+b_n$ for all $n$. Problem \eqref{eq:rpcfnreg} with $\Delta=\Delta_\tau$ can be reformulated as
		\begin{subequations}\label{eq:kvrl}
			\begin{align}
			\min~&~\frac{1}{S}\bd^\top\bone + \mu \bgamma^\top\bone\label{eq:kvrl:cost}\\
			\mathrm{over}~&~\{\bw_n\}_{n=1}^N, \bb, \bd\geq \bzero, \bgamma \label{eq:kvrl:vars}\\
			\mathrm{s.to}~&~-\bar{\bq}_n^g\leq \bZ_n^\top\bw_n +b_n\bone \leq \bar{\bq}_n^g, &&\forall n\label{eq:kvrl:app}\\
			~&~\|\bw_n\|_2\leq \gamma_n,&&\forall n \label{eq:kvrl:norms}\\
			~&~ \|\bX\bq_s^g+\by_s\|_2\leq d_s+\tau,&&\forall s. \label{eq:kvrl:svm}
			\end{align}
		\end{subequations}
		Express voltage deviations at $s$ in terms of $\bw_n$'s and $\bb$
		\begin{equation*}
		\bX\bq_s^g+\by_s = \sum_{n=1}^N \bx_n\bz_{n,s}^\top \bw_n + \sum_{n=1}^N b_n\bx_n + \by_s.
		\end{equation*}
		Let us next introduce the Lagrange multipliers~\cite{SOCP_Lobo}:
		\begin{itemize}
			\item $\ublambda_n\geq \bzero$ and $\oblambda_n\geq \bzero$ corresponding to the linear inequalities in \eqref{eq:kvrl:app} for all $n$;
			\item $(\bu_n,\rho_n)$ related to constraint \eqref{eq:kvrl:norms} for all $n$; and
			\item $(\bmu_s,\sigma_s)$ related to constraint \eqref{eq:kvrl:svm} for all $s$.
		\end{itemize}
		Collect multipliers in {$\bM:=[\bmu_1~\cdots~\bmu_S]\in\mathbb{R}^{N\times S}$}, and vectors $\brho:=[\rho_1~\cdots~\rho_N]^\top$ and  $\bsigma:=[\sigma_1~\cdots~\sigma_S]^\top$. After some algebra, the Lagrangian of \eqref{eq:kvrl} can be written as
		\begin{align}\label{eq:Lagrangian}
		L&=\bd^\top\left(\frac{1}{S}\mathbf{1}-\bsigma\right) + \bgamma^\top\left(\mu\mathbf{1}-\brho\right)\nonumber\\
		&~+ \sum_{n=1}^N \bw_n^\top \left[\bZ_n\left(\oblambda_n - \ublambda_n -\bM^\top \bx_n\right) -\bu_n\right]\nonumber\\
		&~+\sum_{n=1}^N b_n\left[\left(\oblambda_n - \ublambda_n - \bM^\top\bx_n\right)^\top\bone \right]\nonumber\\
		&~-\sum_{n=1}^N \left(\oblambda_n + \ublambda_n\right)^\top  \bar{\bq}_n^g-\sum_{s=1}^S\bmu_s^\top \by_s-\tau\bsigma^\top \bone.
		\end{align}
		Minimizing $L$ over the primal variables provides
		\begin{subequations}\label{eq:lo}
			\begin{align}
			\bsigma&\leq \bone&\label{eq:lo:sigma}\\
			\brho&=\mu\bone&\label{eq:lo:rho}\\
			\bu_n&=\bZ_n \left(\oblambda_n - \ublambda_n - \bM^\top\bx_n\right),&\forall n\label{eq:lo:un}\\
			\left(\oblambda_n - \ublambda_n\right)^\top\bone &=\bx_n^\top\bM\bone,\quad \forall n.\label{eq:lo:lambda}
			\end{align}
		\end{subequations}
		
		From \eqref{eq:lo}, the dual of \eqref{eq:kvrl} becomes the SOCP problem
		\begin{subequations}\label{eq:kvrld}
			\begin{align}
			\max~&-\sum_{n=1}^N \left(\oblambda_n + \ublambda_n\right)^\top  \bar{\bq}_n^g-\sum_{s=1}^S\bmu_s^\top \by_s-\tau\bsigma^\top \bone\label{eq:kvrld:cost}\\
			\mathrm{over}~&~\{\ublambda_n,\oblambda_n\}_{n=1}^N, \{\bmu_s,\sigma_s\}_{s=1}^S \label{eq:kvrld:vars}\\
			\mathrm{s.to}~&~\ublambda_n\geq \bzero,~\oblambda_n\geq \bzero,~ \eqref{eq:lo:lambda},\quad\forall n\label{eq:kvrld:app}\\
			~&~\|\bZ_n \left(\oblambda_n - \ublambda_n - \bM^\top\bx_n\right)\|_2\leq \mu,\quad\forall n \label{eq:kvrld:norms}\\
			~&~\|\bmu_s\|_2\leq \sigma_s\leq 1,\quad\forall s. \label{eq:kvrld:svm}
			\end{align}
		\end{subequations}
		It is not hard to check that \eqref{eq:kvrl} and \eqref{eq:kvrld} are strictly feasible, so strong duality holds and both problems are solvable. The optimal primal and dual variables satisfy complementary slackness SOCPs; see~\cite[Sec.~4.1]{SOCP_Lobo}. For constraints \eqref{eq:kvrl:norms} and \eqref{eq:kvrld:norms}, these conditions identify three cases:
		\renewcommand{\labelenumi}{\emph{c\arabic{enumi})}}
		\begin{enumerate}
			\item If $\|\bw_n\|_2<\gamma_n$, then $\|\bu_n\|_2=\rho_n=0$;
			\item If $\|\bu_n\|_2<\rho_n$, then $\|\bw_n\|_2=\gamma_n=0$; or
			\item If $\|\bw_n\|_2=\gamma_n$ and $\|\bu_n\|_2=\rho_n$, then $\gamma_n\bu_n=-\rho_n\bw_n$.
		\end{enumerate}
		
		Recall that $\rho_n=\mu>0$ from \eqref{eq:lo:rho}. Moreover, it is not hard to see that $\|\bw_n\|_2=\gamma_n$ at the optimum of \eqref{eq:kvrl}. Then, case \emph{c1)} cannot occur. The other two cases entail that $\bw_n=\alpha_n\bu_n$ for some $\alpha_n\leq 0$. Substituting $\bu_n$ from \eqref{eq:lo:un}, and evaluating rule $n$ at the tested scenarios gives
		\begin{align*}
		\bq_n^g&=\bZ_n^\top\bw_n+b_n\bone\\
		&=\alpha_n\bZ_n^\top\bZ_n \left(\oblambda_n - \ublambda_n - \bM^\top\bx_n\right) +b_n\bone\\
		&=\bK_n\ba_n+b_n\bone.
		\end{align*}
		Here we identify $\bK_n=\bZ_n^\top\bZ_n$ and the coefficients in \eqref{eq:rule} as 
		\begin{equation}\label{eq:coeff}
		\ba_n:=\alpha_n \left(\oblambda_n - \ublambda_n - \bM^\top\bx_n\right).
		\end{equation}
		
		Focus now on the complementary slackness for \eqref{eq:kvrl:svm} and \eqref{eq:kvrld:svm}. The equivalent to condition \emph{c1)} reads now as:
		\renewcommand{\labelenumi}{\emph{c\arabic{enumi}')}}
		\begin{enumerate}
			\item If $\|\bX \bq_s^g +\by_s\|_2<d_s+\tau$, then $\|\bmu_s\|_2=\sigma_s=0$.
		\end{enumerate}
		
		Suppose the optimal primal variables satisfy $\|\bX \bq_s^g +\by_s\|_2< \tau$. Then $d_s=0$ follows from \eqref{eq:kvrl}, and \emph{c1')} gives $\|\bmu_s\|_2=\sigma_s=0$. The $s$-th entry of $\ba_n$ in \eqref{eq:coeff} is
		\begin{equation}\label{eq:ant}		a_{n,s}=\alpha_n\left(\overline{\lambda}_{n,s}-\underline{\lambda}_{n,s}-\bmu_s^\top\bx_n\right).
		\end{equation}
		Complementary slackness for \eqref{eq:kvrl:app} implies that $\overline{\lambda}_{n,s} = \underline{\lambda}_{n,s} = 0$ if $|q_{n,s}^g|< \bar{q}_{n,s}^g$ at the optimal, thus proving the claim for linear rules. The result in \eqref{eq:ant} holds for nonlinear rules too. The analysis carries over upon matching the length of $\bw_n$ with the length of $\bphi(\bz_n)$, and substituting $\bZ_n^\top\bZ_n$ by $\bK_n$.
	\end{IEEEproof}		

	\begin{IEEEproof}[Proof of Proposition~\ref{pro:sparsity2}] 
Rewrite \eqref{eq:rpcfnreg} for $\Delta=\Delta_\epsilon$ as
		\begin{subequations}\label{eq:kvrl1}
	   \begin{align}
     	\min~&~\frac{1}{S} \sum_{s=1}^S\bd_s^\top\bone + \mu \bgamma^\top\bone\label{eq:kvrl1:cost}\\
    	\mathrm{over}~&~\{\bw_n\}_{n=1}^N, \bb, \bd\geq \bzero, \bgamma \label{eq:kvrl1:vars}\\
	    \mathrm{s.to}~&~-\bar{\bq}_n^g\leq \bZ_n^\top\bw_n +b_n\bone \leq \bar{\bq}_n^g, &&\forall n\label{eq:kvrl1:app}\\
	   ~&~\|\bw_n\|_2\leq \gamma_n,&&\forall n \label{eq:kvrl1:norms}\\
	   ~&~ -\bd_s-\epsilon \bone \leq \bX\bq_s^g+\by_s \leq \bd_s+\epsilon \bone,&&\forall s. \label{eq:kvrl1:svm}
	  \end{align}
       \end{subequations}
The Lagrangian multipliers of~\eqref{eq:kvrl1} are similar to shose of \eqref{eq:kvrl}, except for {$(\bmu_s,\sigma_s)$} being replaced by ($\underline{\bmu}_s,\overline{\bmu}_s$) and collected in $\underline{\bM}:= [\underline{\bmu}_1~\cdots~\underline{\bmu}_S]$ and $\overline{\bM}:= [\overline{\bmu}_1~\cdots~\overline{\bmu}_S]$.
	Minimizing the Lagrangian of \eqref{eq:kvrl1} over the primal variables yields
	\begin{align*}
		\bu_n= \bZ_n \left[\oblambda_n-\ublambda_n+\left( \overline{\bM}-\underline{\bM} \right) ^{\top} \bx_n \right], \quad \forall n.
	\end{align*}
Similar to Prop.~\ref{pro:sparsity}, the $s$-th entry of $\ba_n$ becomes 
	\begin{align*}
		a_{n,s}=\overline{\lambda}_{n,s}-\underline{\lambda}_{n,s} + (\obmu_s-\ubmu_s)^\top \bx_n.
	\end{align*}
	If the optimal primal variables satisfy {$\|\bX \bq_s +\by_s\|_\infty > \epsilon$, then $\bd_s\neq \bzero$} and accordingly, complementary slackness for~\eqref{eq:kvrl1:svm} implies that $\overline{\bmu}_s\neq \bzero$ or $\underline{\bmu}_s\neq \bzero$. 
\end{IEEEproof}

\bibliographystyle{IEEEtran}
\bibliography{myabrv,power}
\end{document}